\documentclass[a4paper]{amsart}

\usepackage[all]{xy}
\usepackage{amsmath,amssymb,mathpazo}
\SelectTips{cm}{10}

\DeclareMathOperator{\Mod}{Mod}

\renewcommand {\S} {\mathbf S}
\newcommand {\C}{\mathcal C}
\newcommand {\A}{\mathcal A}

\newcommand{\I}{\mathcal I}
\newcommand{\J}{\mathcal J}

\newcommand {\F}{\mathbf F}

\newcommand{\Z}{\mathbf Z}

\newcommand{\op}{{\operatorname{op}}}
\newcommand{\map}{\operatorname{map}}

\newcommand{\id}{\operatorname{id}}

\newcommand {\sm}{\wedge}

\expandafter\def\csname c@figure\endcsname{\csname c@equation\endcsname}
\numberwithin{equation}{section}
\numberwithin{figure}{section}

\newtheorem{lemma}[equation]{Lemma}
\newtheorem{thm} [equation]{Theorem}
\newtheorem{corollary} [equation]{Corollary}
\newtheorem{prop} [equation]{Proposition}

\theoremstyle{definition}
\newtheorem{example}[equation]{Example}

\newtheorem*{defn}{Definition} 
\newtheorem{remark}[equation]{Remark}
\newtheorem*{notation}{Notation}

\DeclareMathOperator{\Cotor}{Cotor}
\DeclareMathOperator{\Tor}{Tor}

\DeclareMathOperator{\coker}{coker}

\DeclareMathOperator{\Hom}{Hom}

\newcommand{\holim}{\operatornamewithlimits{holim}}
\newcommand{\hocolim}{\operatornamewithlimits{hocolim}}
\newcommand{\colim}{\operatornamewithlimits{colim}}

\newcommand{\RKan}{\operatornamewithlimits{RKan}}

\DeclareMathOperator{\Tot}{Tot}

\def\smashop#1_#2{%
\displaystyle{#1_{%
\hbox to 0pt{\hss$\scriptstyle{#2}$\hss}}\;}}
\makeatother
\newdir { >} {{}*!/-5pt/@{>}}

\DeclareMathOperator{\Top}{Top}

\DeclareMathOperator{\im}{im}

\DeclareMathOperator{\Ind}{Ind}

\DeclareMathOperator{\ev}{ev}
\DeclareMathOperator{\odd}{odd}
\newcommand{\N}{\mathbf{N}}

\newcommand{\II}{\mathbf{I}}

\DeclareMathOperator{\diag}{diag}

\makeatletter
\DeclareRobustCommand{\doubleindex}[1]{%
    {
     \edef\resetfontdimens{\noexpand%
         \fontdimen16\textfont2=\the\fontdimen16\textfont2
         \fontdimen17\textfont2=\the\fontdimen17\textfont2\relax}%
     \fontdimen16\textfont2=2.7pt \fontdimen17\textfont2=2.7pt
     #1
     \resetfontdimens}}
\makeatother

\newcommand{\exterior}{\textstyle\bigwedge}
\newcommand{\KK}{\mathcal K}

\newcommand{\QQ}{\mathcal Q}
\newcommand{\EMH}{\underline H}
\DeclareMathOperator{\Pro}{Pro}
\DeclareMathOperator{\Sp}{Sp}
\DeclareMathOperator{\fib}{fib}
\newcommand{\T}{\mathcal T}
\newcommand{\Fin}{\mathcal F}

\title[Convergence of the Eilenberg-Moore spectral sequence]{Convergence of the Eilenberg-Moore spectral sequence for generalized cohomology theories}
\author{Tilman Bauer}
\address{Mathematisches Institut der Universit\"at M\"unster\\
Einsteinstr. 62\\
48149 M\"unster, Germany}
\email{tbauer@math.uni-muenster.de}

\date{March 26, 2008}
\begin{document}

\begin{abstract}
We prove that the Morava-$K$-theory-based Eilenberg-Moore spectral sequence has good convergence properties whenever the base space is a $p$-local finite Postnikov system with vanishing $(n+1)$st homotopy group. 
\end{abstract}

\maketitle

\section{Introduction}
Let $F \to E \to B$ be a fibration of topological spaces. There are three classical spectral sequences allowing one to derive the singular homology of any one of these three spaces from the singular homology of the other two. If $H_*B$ and $H_*F$ are known (the latter as a module over $\pi_1B$), the Serre spectral sequence $H_*(B;H_*(F)) \Longrightarrow H_*(E)$ is a first quadrant spectral sequence which always converges. Slightly less known is the bar spectral sequence (also known as the Rothen\-berg--Steenrod spectral sequence \cite{rothenberg-steenrod}), which allows one to compute $H_*(B)$ whenever $F \to E \to B$ is a principal fibration. In this case, $H_*(F)$ is a ring, $H_*(E)$ is a module over $H_*(F)$, and the natural filtration of the bar construction on $F$ gives a spectral sequence $\Tor^{H_*(F;\F_p)}(H_*(E;\F_p),\F_p) \Longrightarrow H_*(B;\F_p)$ which again is a first-quadrant spectral sequence with good convergence properties. This spectral sequence exists and converges with any coefficients (not just field coefficient), but for the description of the $E^2$-term as a $\Tor$ group one needs a K\"unneth isomorphism. In fact, one can replace $H_*(-;\F_p)$ by any generalized homology theory having K\"unneth isomorphisms and still obtain a strongly convergent right half plane spectral sequence.

This paper is about the dual of the bar spectral sequence, the Eilenberg-Moore spectral sequence (EMSS)
\[
\Cotor^{K_*(B)}_{**}(K_*(E),K_*) \Longrightarrow K_*(F).
\]
Historically predating the bar spectral sequence, the EMSS is actually much harder to understand because it is a second-quadrant (or left half-plane, for nonconnective theories $K$) spectral sequence; in general, it does not converge to its target in any sense. Our main result is:

\begin{thm} \label{mainthm}
Let $p$ be an odd prime, $K(n)$ the $n$th Morava $K$-theory at the prime $p$, and $E_1 \to B \leftarrow E_2$ be a diagram of spaces such that $\pi_*B$ is a finite (graded) $p$-group and $\pi_{n+1}(B)=0$. Let $F = \holim(E_1 \to B \leftarrow E_2)$. Then the $K(n)$-based Eilenberg-Moore spectral sequence,
\[
E^2_{**} = \Cotor^{K(n)_*(B)}_{**}(K(n)_*(E_1),K(n)_*(E_2)) \Longrightarrow K(n)_*F
\]
Ind-converges for any $E_1$, $E_2$. In particular, if $E_1$, $E_2$ are of the homotopy type of finite CW-complexes, the above spectral sequence converges pro-constantly to $K(n)_*F$, which has to be a finite $K(n)_*$-module.
\end{thm}

Before discussing this result, a few words about the history of the problem are in order. The case of $K_*=H_*(-;\F_p)$ has been studied extensively \cite{eilenberg-moore:homology-and-fibrations,smith:emss,dwyer:strongconvergence,dwyer:exoticconvergence}, and the convergence issues arising here are the same as for any \emph{connective} theory $K$. Roughly, the question of convergence only depends on $\pi_1(B)$ and its action on the homology of $F$.

Now let $K$ be nonconnective and possessing K\"unneth isomorphisms. Thus $K$ is one of Morava's extraordinary $K$-theories $K(n)$ or an extension of it, which is the case I am interested in in this work. The question of convergence becomes much more intricate; in particular, nonconvergence can occur even for simply connected base spaces. As an example, consider the path-loop fibration $K(\Z/2,1) \to * \to K(\Z/2,2)$ and $K=K(1) = KU/2$, mod-2 ordinary $K$-theory. In this case $K(1)_*(K(\Z/2,2)) = K(1)_*$, but $K(1)_*(K(\Z/2,1))$ is nontrivial, so that there is no chance for the (trivial) EMSS to converge. Theorem~\ref{mainthm} says that this nonvanishing of $\pi_{n+1}$ is in fact the only obstruction to convergence if the base space has totally finite homotopy groups which are $p$-groups. Note that the condition that the homotopy groups are $p$-groups is not too restrictive because we can always replace the fibration under consideration by its $\Z/p$-localization, and the EMSS will never know the difference (although the fiber might change drastically).

Previous work on the $K(n)$-based EMSS includes work by Tamaki \cite{tamaki:emss}, where he shows convergence when the base space is of the form $\Omega^{n-1}\Sigma^n X$ (he mistakenly claims strong convergence), and work by Jeanneret and Osse \cite{jeanneret-osse:emss}, where the authors show convergence whenever the base space has certain homological global finiteness properties, for example, if the base space is the classifying space of a polynomial $p$-compact group.

The term ``Ind-convergence'' in Theorem~\ref{mainthm} requires explanation. We first recall the classical notion of \emph{pro-convergence} (called strong convergence in \cite{bousfield:homology-cosimplicial,shipley:homology-cosimplicial}, but different from Cartan-Eilenberg's and Boardman's notion of strong convergence \cite{ce:homalg,boardman:conditionally}). Associated to $F = \holim(E_1 \to B \leftarrow E_2)$, there is a tower of $K(n)$-module spectra $T^\bullet(E_1,E_2) = \Tot^\bullet K(n)[C_B(E_1,E_2)]$ coming from the two-sided cobar construction of $E_1$ and $E_2$ over $B$, and a map $K(n)[F] \to T^\bullet(E_1,E_2)$. The spectral sequence always converges conditionally to $\pi_*(\holim T^\bullet)$, which may be different from $K(n)_*(F)$. We say that the spectral sequence is pro-constantly convergent if $K(n)_*F \to \pi_*T^\bullet$ is a pro-isomorphism from the constant object $K(n)_*F$ to this natural target. This in particular implies that the spectral sequence converges in a very strong sense, namely, only finitely many differentials live at any bidegree, and in $E^\infty$ the filtration is finite in every total degree.

We consider the EMSS not as one spectral sequence, but as a whole directed system of spectral sequences, one for each pair of finite sub-CW-complexes of $E_1$ and $E_2$. Similarly, the target $K(n)_*F$ can be thought of as the directed system of $K(n)_*F'$ where $F'$ runs through all finite sub-CW-complexes of $F$. We call the spectral sequence Ind-convergent if
the comparison map $K(n)_*F \to \pi_*T^\bullet(E_1,E_2)$ is an isomorphism in the category of ind-pro-abelian groups. In particular, if $E_1$ and $E_2$ are finite CW-complexes, then $T^\bullet(E_1,E_2)$ is ind-constant, thus $K(n)_*F$ also has to be ind-constant, which can only happen if $K(n)_*(F)$ is finite. Furthermore, in this case, since $K(n)_*F$ is pro-constant, $T^\bullet$ also has to be isomorphic to a pro-constant tower, and we get the specialization mentioned in the theorem.

We also mention that Ind-convergence is a good enough notion to allow for a comparison theorem of spectral sequences:
\begin{corollary}
Let $(E_1' \to B' \leftarrow E_2') \to (E_1 \to B \leftarrow E_2)$ be a map of diagrams such that the associated spectral sequences $(E^*)'$, $E^*$ converge Ind-constantly to $K_*(F')$ and $K_*(F)$, respectively. Then there is an induced ind-map $(E^*)' \to E^*$, and if for any $s$, $(E^s)' \to E^s$ is an ind-isomorphism, then $K_*(F) \cong K_*(F')$.
\end{corollary}
We get an induced map because every finite sub-CW-complex of $E_i'$ maps to another finite sub-CW-complex of $E_i$. The result follows from the fact that we get an ind-isomorphism $K_*(F) \to K_*(F')$, and the ordinary group $K_*(F)$ is simply the colimit of the ind-group $K_*(F)$.

The heuristic reason for introducing this slightly unwieldy notion of Ind-con\-ver\-gence is that the natural target of the spectral sequence, $\holim T^\bullet$, does not commute with infinite colimits or even infinite coproducts because it is an inverse limit. However, the fiber of a colimit of total spaces over a fixed base space is just the colimit of the total spaces of the individual fibrations, and thus we have to ``train'' the spectral sequence to commute with colimits. This is achieved by passing to the ind-category.

Theorem~\ref{mainthm} is proved by means of a series of results of possibly independent interest. We collect the main steps here. Let $K=K(n)$ be the $n$th Morava $K$-theory.

\begin{thm} \label{hopfringconvergence}
For odd $p$, the $K(n)$-based Eilenberg-Moore spectral sequence Ind-converges for the path-loop fibration on $K(\Z/p,m)$ whenever $m \neq n+1$.
\end{thm}

This result is proved in Section~\ref{cohopf} by a complete computation of the spectral sequence, aided by the computations of $K(n)_*(K(\Z/p,m))$ for all $m$, $n$ in \cite{ravenel-wilson:morava}. It is likely that the result also holds for $p=2$, but extra care is needed due the non-commutativity of $K(n)$ in this case.

To pass from a contractible total space to more general cases, in Section~\ref{indeptotal} we prove:
\begin{thm} \label{universalconvergence}
Let $B$ be a space such that $* \to B \leftarrow *$ has an Ind-convergent EMSS. Then so has $E_1 \to B \leftarrow E_2$ for any $E_1$, $E_2$.
\end{thm}

We call a space $B$ with this property \emph{Ind-convergent}.

Finally, in Section~\ref{transitivitysection} we show how to pass to more complicated base spaces than just $K(\Z/p,m)$:

\begin{thm} \label{fibrationconvergence}
Let $F \to Y \to X$ be a fibration with $F$ and $X$ Ind-convergent. Then so is $Y$.
\end{thm}

\begin{proof}[Proof of Theorem~\ref{mainthm}]
By Theorem~\ref{universalconvergence}, it suffices to consider the case $E_1=E_2=*$. Since $B$  has finite homotopy groups which are $p$-groups, it has a finite Postnikov decomposition
\[
\xymatrix{
B = B_k \ar[d]\\
\vdots\ar[d]\\
B_1 \ar[d] \ar[r] & K(\Z/p,n_3)\\
B_0 = K(\Z/p,n_1) \ar[r] & K(\Z/p,n_2),
}
\]
where none of the $n_i$ equal $n+1$. By induction, Theorem~\ref{fibrationconvergence} and Theorem~\ref{hopfringconvergence}, we conclude that $B$ is Ind-convergent.
\end{proof}

Section~\ref{emssintro} is an exposition of the construction of the Eilenberg-Moore spectral sequence in the generality we need, and Section~\ref{indprosection} deals with the structure of ind-pro-objects.

\section{The generalized Bousfield and Eilenberg-Moore spectral sequences} \label{emssintro}

Associated with any cosimplicial space $C^\bullet$ there is a tower of spaces,  called the \emph{Tot-tower},
\[
\{\Tot^s C^\bullet\}_{s\geq 0}, \quad \Tot^s C^\bullet = \map(\Delta^\bullet_{\leq s},C^\bullet),
\]
where $\Delta^\bullet$ denotes the cosimplicial space whose $s$th space $\Delta^s = |\Delta[s]|$ is the standard $s$-simplex, $\Delta^\bullet_{\leq s}$ denotes the $s$-skeleton of $\Delta^\bullet$, and ``$\map$'' is the mapping space functor. This tower is the target of a canonical map from $\Tot C^\bullet =  \holim_{s} \Tot^s C^\bullet = \map(\Delta^\bullet,C^\bullet)$. A similar construction can be made for cosimplicial spectra; these are connected by homotopy equivalences
\[
\Tot^s \Omega^\infty C^\bullet \simeq \Omega^\infty \Tot^s C^\bullet,
\]
for cosimplicial spectra $C^\bullet$, but in general,
\[
\Tot^s \Sigma^\infty C^\bullet \not\simeq \Sigma^\infty \Tot^s C^\bullet.
\]
More generally, let $K$ be any spectrum, and define $K[X]$ to be the $K$-module spectrum $K \sm \Sigma^\infty(X_+)$ on the space $X$. There is a natural map of towers of spectra
\[
\{ K[\Tot^s C^\bullet] \}_s\xrightarrow{\Phi} \{\Tot^s(K[C^\bullet])\}_s
\]
which is rarely a homotopy equivalence.

For both towers there is a spectral sequence abutting to their respective homotopy inverse limits; it is, however, only the one on the right hand side whose $E^1$- and $E^2$-terms have a convenient formulation. For we have a cofibration sequence of spectra
\begin{equation}\label{normalizedcofseq}
\Sigma^{-s}N^s(K[C^\bullet]) \to \Tot^s(K[C^\bullet]) \to \Tot^{s-1}(K[C^\bullet])
\end{equation}
where the normalization $N^s$ is the fiber of $C^s \to M^{s-1} (K[C^\bullet])$, the latter being the cosimplicial matching space \cite[Chapters VII.4, VIII.1]{goerss-jardine}. This yields $E^2_{s,t} = \pi^s K_t(C^\bullet)$. This spectral sequence is the \emph{$K$-based Bousfield spectral sequence} or homology spectral sequence of a cosimplicial space \cite{bousfield:homology-cosimplicial}.

The spectral sequence belonging to the tower in the source of $\Phi$, however, has a more accessible target, namely $\pi_* \holim_s K[\Tot^s C^\bullet]$. There is a map of towers
\begin{equation}\label{Pmap}
\{K[\Tot C^\bullet]\} \xrightarrow{P} \{ K[\Tot^s C^\bullet] \}_s,
\end{equation}
where the left tower is constant, and thus we get a comparison map from $K_*\Tot C^\bullet$ to the target of either spectral sequence.

Bousfield \cite{bousfield:homology-cosimplicial} studied when the maps $P$ and $\Phi$ are pro-isomorphisms in the case of $K=H\F_p$, which would imply that the associated spectral sequences all converge to $K_*(\Tot C^\bullet)$. Various criteria were given for convergence \cite[Theorems 3.2, 3.4, 3.6]{bousfield:homology-cosimplicial}, which generalize to the case of connective homology theories $K$. However, convergence for periodic theories remained an intricate problem.

\subsection{Forms of convergence}

Let us examine the notion of convergence in the Bousfield spectral sequence more closely. We define a decreasing filtration $F^\bullet K_*(\Tot C^\bullet)$ by
\[
F^sK_*(\Tot C^\bullet) = \ker\left( K_*(\Tot C^\bullet) \xrightarrow{\Phi_* \circ P_*} \pi_* \Tot^s (K[C^\bullet])\right).
\]

Recall from \cite{ce:homalg,boardman:conditionally} that the spectral sequence is called 
\emph{strongly convergent} to $K_*(\Tot C^\bullet)$ if two conditions are satisfied:
\begin{enumerate}
\item The natural map $F^s / F^{s+1} K_*(\Tot C^\bullet) \to E_\infty^{s,*}$ is an isomorphism and
\item the filtration $F^\bullet K_*(\Tot C^\bullet)$ is complete Hausdorff, i.~e.
\[
\lim F^s = {\lim}^1 F^s = 0.
\]
\end{enumerate}
It is called \emph{completely convergent} if it is strongly convergent and additionally,
\[
{\lim_s}^1 \pi_* \Tot^s (K[C^\bullet]) = 0.
\]

\begin{remark}\label{complconvremark}
Assuming strong convergence, complete convergence is equivalent to $K[\Tot C^\bullet] \simeq \Tot K[C^\bullet]$: In fact, strong convergence means that
\[
K_*(\Tot C^\bullet) \cong \lim_s \pi_*\Tot^s(K[C^\bullet]);
\]
the Milnor exact sequence
\[
0 \to {\lim_s}^1 \pi_{t+1} \Tot^sK[C^\bullet] \to \pi_t \Tot(K[C^\bullet]) \to \lim_s \pi_t \Tot^s(K[C^\bullet] \to 0
\]
thus shows that the $\lim^1$-term vanishes if and only if $K[-]$ commutes with $\Tot$ for $C^\bullet$.
\end{remark}

The litmus test for the usability of any notion of convergence is whether it implies the spectral sequence comparison theorem. For strong convergence, this was proved in \cite[Theorem 5.3]{boardman:conditionally}:
\begin{thm}[Boardman]
Let $f\colon C^\bullet \to D^\bullet$ be a map of cosimplicial spaces such that the Bousfield spectral sequences for $C^\bullet$ and $D^\bullet$ converge strongly. If $f$ induces an isomorphism on any $E^s$-term ($1 \leq s \leq \infty$) then $f$ also induces an isomorphism $K_*(\Tot C^\bullet) \to K_*(\Tot D^\bullet)$.
\end{thm}

There is a different and stronger version of the term strong convergence in the context of the Bousfield or Eilenberg-Moore spectral sequence. In \cite{bousfield:homology-cosimplicial,shipley:homology-cosimplicial}, the Bousfield spectral sequence associated to a homology theory $K$ and a cosimplicial space $C^\bullet$ is called strongly convergent if the tower map
\[
\left\{K_t(\Tot C^\bullet)\right\}_s \xrightarrow{\Phi_* \circ P_*} \left\{\pi_t \Tot^s (K[C^\bullet])\right\}_s
\]
is a pro-isomorphism for each $t$, where the tower on the left hand side is constant. Explicitly, this means that for every $t \in \Z$ and $s \geq 0$, there is an $N(s,t) \in \N$ and a map in the following diagram, making both triangles commute:
\[
\xymatrix{
K_t(\Tot C^\bullet) \ar@{=}[d] \ar[r]^-{\Phi_* \circ P_*} & \pi_t\Tot^{s+N(s,t)}K[C^\bullet] \ar@{-->}[dl] \ar[d]\\
K_t(\Tot C^\bullet) \ar[r]^{\Phi_* \circ P_*} & \pi_t\Tot^s K[C^\bullet]
}
\]

We will call this kind of convergence \emph{pro-constant convergence}. If the function $N(s,t)$ can be chosen to be constant (say $\equiv N$), we call the tower map an \emph{$N$-isomorphism} and the spectral sequence \emph{$N$-convergent}. In our applications, $K$ is a periodic homology theory, which means that if the spectral sequence is pro-constantly convergent, we can choose $N(s,t)$ to be independent of $t$; viz, take
\[
N(s) = \max\{N(s,t) \mid 0 \leq t < \text{period of $K$}\}.
\]

Recall the following well-known lemma, which is a generalization of \cite[Lemma~3.5]{bousfield:homology-cosimplicial}:

\begin{lemma} \label{spseqconvergence}
The Bousfield spectral sequence $E^*_{**}$ associated to a tower of spectra $C^\bullet$ is pro-constantly convergent if and only if
\begin{enumerate}
\item For each $s,\; t$, there is an $N=N(s,t)>0$ such that $E^N_{s,t}=E^\infty_{s,t}$ in the spectral sequence and
\item For each $k$ there is an $N(k)$ such that $E^\infty_{s,s+k}=0$ for $s\leq N(k)$.
\end{enumerate}
Moreover, the spectral sequence is $N$-convergent if and only if $N(s,t)$ and $N(k)$ above can be chosen to be constant with value $N$.
\end{lemma}

\begin{lemma} \label{proconstimpliescomplete}
Pro-constant convergence implies complete convergence.
\end{lemma}
\begin{proof}
Let $Y \to X^\bullet$ be a tower of spectra under $Y$ such that the associated spectral sequence is pro-constantly convergent. By Lemma \ref{spseqconvergence}, this implies that for each $s,t$, there is an $n$ such that $E_n^{s,t} = E_\infty^{s,t}$. Thus the \emph{derived $E_\infty$ term} $RE_\infty = \lim_r^1 Z_r$ is zero, and by \cite{boardman:conditionally}, the spectral sequence is completely convergent.
\end{proof}

\begin{remark}\label{Nconvremark}
The Bousfield spectral sequence associated to a tower of spectra $C^\bullet$ is pro-constantly (resp. $N$-) convergent if and only if
\begin{enumerate}
\item $\Phi\colon K[\Tot C^\bullet] \to \Tot K[C^\bullet]$ is a homotopy equivalence; and
\item The tower $\pi_t\Tot K[C^\bullet]$ is pro- (resp. $N$-) constant for each $t$.
\end{enumerate}
The reason for this, as in Remark~\ref{complconvremark}, is that for pro-constant towers, the derived functor of the inverse limit is trivial.
\end{remark}

While complete convergence is a perfectly fine property to ensure that the spectral sequence determines its target uniquely up to filtration, it has the technical disadvantage that the levelwise cofiber of two completely convergent towers need not be completely convergent, making it all but unusable for inductive arguments. On the other hand, the tower Five Lemma \cite{yellowmonster} implies that pro-constant convergence is preserved under taking levelwise cofibers. In general, there is no Five Lemma for $N$-isomorphisms.

\subsection{The Eilenberg-Moore spectral sequence as a special case of the Bousfield spectral sequence}
Let
\[
\xymatrix{
F \ar[r] \ar[d] &  E_1 \ar[d]^{p_1} \\ E_2 \ar[r]^{p_2} & B
}
\]
be a pullback diagram of spaces, where at least one of $E_1 \to B$ and $E_2 \to B$ is a fibration. Denote by $C^\bullet=C^\bullet_B(E_1,E_2)$ the cobar construction for $E_1 \to B \leftarrow E_2$; that is, $C^s = E_1 \times B^s \times E_2$ with the usual coboundary and codegeneracy maps. We abbreviate $C_B(E)$ for $C_B(E,*)$. The Bousfield spectral sequence associated to $C^\bullet$ is called the Eilenberg-Moore spectral sequence. The convergence question in this special case becomes a little bit simpler, but is still hard to tackle. In fact, one readily sees that $\Tot^0 C^\bullet = E_1 \times E_2$, whereas $\Tot^i C^\bullet \simeq F$ for all $i \geq 1$. Thus the map $P$ of \eqref{Pmap} is always an isomorphism in positive degrees. In Bousfield's terms, in this case ``pro-convergence'' implies pro-constant convergence.

The description of the $E^1$-term becomes significantly simpler than in the generic Bousfield spectral sequence since one can avoid the computation of the normalization of the cosimplicial spectrum $K[C^\bullet]$ as in \eqref{normalizedcofseq}.

\begin{lemma} \label{totfibers}
The cosimplicial spectrum $K[C^\bullet]$ is \emph{codegeneracy-free}, i.e. it is the right Kan extension of a diagram $Z\colon \Lambda \to \{\text{Spectra}\}$, where $\Lambda$ is the category with objects $[n]$, $n \geq 0$, and injective monotonic maps as morphisms. This diagram is defined by $Z^n = K \sm \left((E_1)_+ \sm B^{\sm n} \sm (E_2)_+\right)$.
\end{lemma}
\begin{proof}
The coface maps in $Z^\bullet$ are defined by
\begin{align*}
d^0(e_1,b_1,\dots,b_n,e_2) &= (e_1,p_1(e_1),b_1,\dots,b_n,e_2),\\
d^{n+1}(e_1,b_1,\dots,b_n,e_2)  &= (e_1,b_1,\dots,b_n,p_2(e_2),e_2)\\
\intertext{and}
d^i(e_1,b_1,\dots,b_n,e_2) &= (e_1,b_1,\dots,b_i,b_i,\dots,b_n,e_2) \text{ for }1 \leq i \leq n.
\end{align*}
Let $I\colon \Lambda \to \Delta$ denote the inclusion functor. Then the right Kan extension over $I$ of any $\Lambda$-diagram of spectra $Z^\bullet$ can explicitly be described as
\[
{\RKan_I Z^\bullet}^n = \bigvee_{[n] \twoheadrightarrow [k]} Z^k
\]
where the wedge runs over all surjections $[n] \to [k]$ (of which there are $n \choose k$).

For every such $\phi\colon [n] \to [k]$, there is a map $\phi^*\colon B^{\sm k} \to (B_+)^{\sm k} = B^n_+$ given on the $i$th coordinate ($i=1,\dots,n$) by
\[
\phi^*_i(b_1,\dots,b_k) = \begin{cases} *; & \phi(i)=\phi(i-1)\\
b_{\phi(i-1)}; & \text{otherwise}\end{cases}
\]
Now
\begin{align*}
K[C^n] &= K[E_1 \times B^n \times E_2] = K \sm \left((E_1)_+ \sm (B_+)^{\sm n} \sm (E_2)_+\right)\\
&\xleftarrow{(\phi^*)} \bigvee_{\phi\colon [n] \twoheadrightarrow [k]} K \sm \left((E_1)_+ \sm B^{\sm k} \sm (E_2)_+\right)
\end{align*}
is an isomorphism compatible with the cosimplicial structure maps.
\end{proof}

This implies that
\[
E^1_{r,s} = K_{r-s}((E_1)_+ \sm B^{\sm n} \sm (E_2)_+).
\]
The $E^2$-term also has a convenient description, at least if $K_*$ is a graded field or, more generally, if $K_*(B \times X) \cong K_*(B) \otimes_{K_*} K_*(X)$ for all $X$. It is given by
\[
E^2_{r,s} = \Cotor^{K_*B}_{r,s}(K_*E_1,K_*E_2).
\]

\subsection{The Eilenberg-Moore spectral sequence for parametrized spectra} \label{paramspectra}

It is crucial for the Eilenberg-Moore spectral sequence that the base space is indeed a space and not, say, a spectrum. That is, we cannot expect a functorial spectral sequence that takes as input a diagram of spectra $E_1 \to B \leftarrow E_2$ and which computes something that in the case of suspension spectra of spaces is the suspension spectrum of the homotopy pullback. The technical reason is that the diagonal on $B$ is needed for the cobar construction. However, for $E_1$ and $E_2$ we only need the coaction $E_i \to E_i \times B$, never the diagonal. In this section, I will set up a convenient category of $K$-module spectra over a space $B$ and show that the Eilenberg-Moore spectral sequence can be generalized for pairs of such $K$-module spectra over $B$.

Let $B$ be a pointed space. The objects of the category $(\Top/B)_*$ of sectioned spaces over $B$ are of the form $(X \xrightarrow{p_X} B,s_X)$, where $s_X\colon B \to X$ is a section of $p_X$. The maps are given by maps $f\colon X \to Y$ over $B$ compatible with the sections. In the category $(\Top/B)_*$ one can define fiberwise homotopical constructions such as cofibers, fibers, suspensions, smash products etc, which we will denote by adding a subscript $B$ to the usual symbol, e.g. $\sm_B$. The category $(\Top/B)_*$ is also complete and cocomplete. For details about these construction, consult e.g. \cite{smith:emss,may-sigurdsson}.

The category $\Sp_B$ is the category of spectra over the \emph{space} $B$. An object is a sequence $X_i \xrightarrow{\curvearrowleft} B$ in $(\Top/B)_*$ together with maps $\Sigma_B X_i \to X_{i+1}$. If $X \in (\Top/B)_*$ is a sectioned space, the fiberwise suspension spectrum $\Sigma^\infty_B X$ is an object in $\Sp_B$. As for spaces, all homotopical constructions such as (co-)fibers, smash products, (co-)limits work in this fiberwise setting. In \cite{may-sigurdsson}, this is treated in a ``brave new'' way, but for our purposes, the naive notions of spectra are enough.

If $X$ is a pointed space (not over $B$), we denote by $X_B$ the object $(X \times B \to B,s)$, where $s(b) = (*,b)$. Similarly, if $E$ is any spectrum, we denote by $E_B$ the spectrum over $B$ whose $n$th space is $(E_n)_B$. This is indeed a spectrum over $B$ by means of the structure maps
\[
\Sigma_B (E_n)_B = (\Sigma E_n)_B \to (E_{n+1})_B.
\]
In fact, the functor $X \to X_B$ is right adjoint to the forgetful functor $U\colon \Top/B \to \Top$ (forgetting the map to $B$). and $E \to E_B$ is right adjoint to the $U\colon \Sp/B \to \Sp$ sending $(E \to B, s)$ to $\{E_n/s(B)\}$, which is an ordinary spectrum.

In particular, for an ordinary ring spectrum $K$, we have the notion of a fiberwise $K$-module spectrum, i.e. a spectrum $X \in \Sp_B$ with a homotopy associative and unital action $K_B \sm_B X \to X $. We denote the category of $K$-bimodule spectra over $B$ by $\Mod_K/B$. There is a fiberwise smash product $\sm_{K,B}\colon \Mod_K/B \times \Mod_K/B \to \Mod_K/B$. 

Now let $K$ be a ring spectrum and $E_1, E_2 \in \Mod_K/B$. We define the cobar construction $C_B(E_1,E_2)$ to be the cosimplicial $K$-module spectrum
\[
C_B(E_1,E_2)^n = \left(U(E_1) \sm_{K} K[B^n] \sm_{K} U(E_2)\right),
\]
where the middle face and degeneracy maps are induced by the diagonal of $B$ and projections, and the adjunction counits
\[
E_i \to (U(E_i))_B
\]
induce the remaining structure maps
\[
U(E_i) \to U(U(E_i)_B) = U(E_i) \sm_K K[B].
\]
If $E_1 \to B \leftarrow E_2$ is a diagram of spaces, we have that
\[
C^n(K[E_1] \vee B, K[E_2] \vee B) \cong K[E_1] \sm_K K[B^n] \cong K[E_2] = K[C_B^n(E_1,E_2)],
\]
thus we recover our original cobar construction for spaces.

\section{Pro-objects and Ind-Pro-objects} \label{indprosection}

In this section, we will study ind-pro-objects in $K$-module spectra and their homotopical properties. We start with a motivating example illustrating our need to introduce ind-structures.

Consider an infinite set of fibrations $E_i \to B$ with the same base space and inclusions $E_i \to E_{i+1}$, and let $E = \bigcup_i E_i$. Then obviously, the fiber of $E \to B$ is the union of the various fibers $F_i$ of $E_i  \to  B$. Thus, taking fibers commutes with filtered colimits. However, if we study the construction of the Eilenberg-Moore spectral sequence, we find that while
\[
\KK^{(n)}(X) = \Tot^nK[C_B(X)]
\]
still commutes with filtered colimits ($\Tot^n$ is a finite limit),
\[
\KK(X) = \holim_n \KK^{(n)}(X) = \Tot K[C_B(X)]
\]
might not because the inverse limit has no reason to commute with colimits. Thus, in a way, $\KK(X)$ is not ``the correct target'' of the spectral sequence. To offset this deficit, we could think of $X$ as the directed system of all compact subobjects of $X$, i.e., all finite sub-CW-complexes of $X$ if we assume $X$ to be a CW-complex. Applying $\KK$ to this system, we obtain a functor that now commutes with all filtered colimits and thus represents a better target for the spectral sequence, which is now really a filtered diagram of spectral sequences. Thus, instead of looking at towers as objects in the pro-category $\Pro-\Mod_K$, we are now looking at ind-objects in this pro-category, that is, objects in $\T = \Ind-\Pro-\Mod_K$.

Let $\C$ be any category. The category $\Pro-\C$ has as objects pairs $(\I,X\colon \I \to \C)$ where $\I$ is a cofiltered small category and $X$ is a functor. The morphisms are given by
\[
\Hom_{\Pro-\C}((\I,X),(\J,Y)) = \lim_\J \colim_\I \Hom_\C(X,Y).
\]
It is useful to think of this as saying that giving a morphism is giving for every $j \in \J$ a map $X(i) \to Y(j)$ for some $i \in \I$, although this ignores the fact that these have to be compatible in some way.

Dually, the category $\Ind-\C$ has as objects pairs $(\I,X)$ as above but with $\I$ a filtered category; the morphisms are given by
\[
\Hom_{\Ind-\C}((\I,X),(\J,Y)) = \lim_\I \colim_\J \Hom_\C(X,Y).
\]

Recall \cite[Appendix]{artin-mazur:etale} that any map $X \to Y$ in $\Ind-\C$ or $\Pro-\C$ can be represented by a level map, that is, there is a filtered category $\I$ (or, without loss of generality, a directed set $\I$), functors $X',\; Y'\colon \I \to \C$, a natural transformation $X' \to Y'$ and isomorphisms $X \cong X'$, $Y \cong Y'$ in $\Ind-\C$ (resp. $\Pro-\C$) such that $X \to Y$ is the composite $X \cong X' \to Y' \cong Y$. Similarly, any finite, loopless diagram $D \to \Ind-\C$ or $D \to \Pro-\C$ is isomorphic to a diagram of levelwise maps.

However, an object $X \in \Ind-\Pro-\C$ is not necessarily isomorphic to a doubly indexed system $X^t_s$, where $t$ runs through an inverse system $T$ and $s$ runs through a direct system $S$. Any such $X$ is isomorphic to a diagram $S \to \Pro-\C$, for a directed set $S$, but since $S$ is not finite, we cannot replace this diagram by a levelwise diagram. Instead, $X$ can always be represented by a functor $\QQ_{S,\alpha} \to \C$, where $S$ is a directed set, $\alpha\colon S^{\op} \to \{\text{inverse sets}\}$ is a functor, and $\QQ_{S,\alpha}$ is a poset with objects pairs $(s,t)$ ($s \in S, t \in \alpha(s)$) and $(s,t)\leq(s',t')$ if $s\leq s'$ and $t \leq \alpha(s < s')(t') \in \alpha(s)$.

\begin{remark}\label{proproremark}
For later reference, we note that pro-pro-objects can be described in a similar way, giving an inverse set $S$, a functor $\alpha\colon S^{\op} \to \{\text{inverse sets}\}$, and a functor $X\colon \QQ_{S,\alpha} \to \C$. In this case, unlike for ind-pro-objects, the poset $\QQ_{S,\alpha}$ is again an inverse set, and thus $X$ can also be interpreted as an object in $\Pro-\C$. We denote this tautological ``reinterpretation functor'' by $D\colon \Pro-\Pro-\C \to \Pro-\C$. (The letter is supposed to remind one of the diagonal of a double tower.)
\end{remark}

We call a natural transformation of functors $X \to Y\colon \QQ_{S,\alpha}$ in $\Ind-\Pro-\C$ a \emph{levelwise map}. The fact that any ind- oder pro-map is isomorphic to a levelwise map now easily generalizes to
\begin{lemma} \label{levelwise}
Every map in $\Ind-\Pro-\C$ is isomorphic to a levelwise map. \qed
\end{lemma}



In our applications, the category $\C$ will be either $\Top/B$, the category of topological spaces over a base space $B$, or $\Mod_K/B$, the homotopy $K$-module spectra over $B$. Recall from Section~\ref{paramspectra} that there are forgetful functors $U\colon \Top/B \to \Top$ (forgetting the map to $B$) and $U\colon \Mod_K/B \to \Mod_K$ (sending $(X \to B, s)$ to $X/s(B)$). If $K$ is a homology theory or $\pi_*$, we write $K_*(X)$ for $K_*(U(X))$.

The category $\Top/B$ carries a model structure, where a map $f$ is a weak equivalence, fibration, or cofibration if the underlying map $U(f)$ in $\Top$ is a weak equivalence, Serre fibration, or Serre cofibration. However, as is pointed out in great detail in \cite[Ch.~6]{may-sigurdsson}, this model structure has bad properties: for example, the fibrant replacement functor does not in general commute with cofibers, even up to weak equivalence. However, there is a model structure, which May and Sigurdsson call the qf-structure, which has the same weak equivalences, but different cofibrations and fibrations \cite[Thm 6.2.5]{may-sigurdsson}. This allows us to equip the category of spectra over $B$ with a good model structure as well \cite[Thm 12.3.10]{may-sigurdsson} whose fibrant objects are the $\Omega$-spectra over $B$. Similarly, the category of $K$-module spectra over $B$ carries a model structure \cite[Thm 14.1.7]{may-sigurdsson}. In these two structures, the weak equivalences are the stable homotopy equivalences after applying $U$.

The reader should be warned that constructing these model structures is a substantial amount of work. However, for our purposes, it is enough to know that there exists a model structure with the right weak equivalences. We can use this model structure as a black box.

\begin{defn}
A map $X \to Y$ in $\Ind-\C$ or in $\Ind-\Pro-\C$ is an \emph{essentially levelwise weak equivalence} if it is isomorphic to a levelwise map which is a weak equivalence on every level.
\end{defn}

It was proven in \cite{isaksen:strictmodel} that a composite of essentially levelwise weak equivalences is again an essentially levelwise weak equivalence (in fact, that the ind-objects in any proper model category carry a model structure where the weak equivalences are the essentially levelwise weak equivalences). 

%

For the following results, we need $K$ to be a field. In this case, we have that
\[
[X,Y]_K \cong \Hom_{K_*}(\pi_*X,\pi_*Y),
\]
where the left hand side means homotopy classes of $K$-module maps.

\begin{lemma} \label{levelwisewelemma}
Let
\begin{equation}\label{htpyisodiagram}
\xymatrix{
X \ar[r]^{f} \ar[d]^p & Y \ar@{.>}[dl] \ar[d]^p\\
X' \ar[r]^{f'} & Y',
}
\end{equation}
be a commutative diagram in $\Mod_K/B$, where the dotted arrow denotes a map $\pi_*Y \to \pi_*X'$ such that the diagram commutes on homotopy groups.
\begin{enumerate}
\item \label{makesurjective} Assume $X'$ is cofibrant and fibrant and $Y$ is fibrant in $\Mod_K$, and that $f'$ is a fibration in $\Mod_K/B$. Then there exists $\tilde X \in \Mod_K/B$ and a commutative diagram
\[
\xymatrix@C=0.5pc@R=0.5pc{
X \ar[rr]^{f} \ar[dd] \ar[dr] && Y \ar[dd]\\
& \tilde{X} \ar[ur] \ar[dl]\\
X' \ar[rr]^{f'} && Y',
}
\]
where $\tilde X \to Y$ is a fibration inducing a surjective map in homotopy.
\item \label{liftsurjective} If we have a map of diagrams $D_1 \to D$ of the form \eqref{htpyisodiagram}, given by $X_1 \to X$, $Y_1 \to Y$, etc., and $\tilde X$ as in \eqref{makesurjective}, there is a completion $\tilde X_1$ of $D_1$ as in \eqref{makesurjective} and a map $\tilde X_1 \to \tilde X$ making everything commute:
\[
\xy
\xymatrix"0"{
X \ar[r]^f \ar[d] & Y \ar[dd]^p\\
\tilde X \ar[ur] \ar[d]\\
X' \ar[r]_{f'} & Y'
}
\POS-(10,-10)
\xymatrix"1"{
X_1 \ar[r]^{f_1} \ar["0"] \ar[d] & Y_1 \ar[dd]^{p_1} \ar["0"]\\
\tilde X_1 \ar["0"] \ar[ur] \ar[d]\\
X_1' \ar[r]^{f_1'} \ar["0"] & Y_1' \ar["0"]\\
}
\endxy
\]
\item \label{makeinjective} Dually, assume $Y$ is cofibrant and fibrant and $X'$ is cofibrant in $\Mod_K$, and that $f$ is a cofibration in $\Mod_K/B$. Then there exists $\tilde Y' \in \Mod_K/B$ and a commutative diagram
\[
\xymatrix@C=0.5pc@R=0.5pc{
X \ar[rr]^{f} \ar[dd] && Y \ar[dd] \ar[dl]\\
& \tilde{Y'} \ar[dr] \\
X' \ar[rr]^{f'} \ar[ur] && Y',
}
\]
where $X' \to \tilde Y'$ is a cofibration inducing an injective map in homotopy.
\item \label{liftinjective} If we have a map of diagrams $D_1 \to D$ as in \eqref{liftinjective} and $\tilde X$ as in \eqref{makeinjective}, there is a completion $\tilde X_1$ of $D_1$ as in \eqref{makeinjective} and a map $\tilde X_1 \to \tilde X$ making everything commute:
\[
\xy
\xymatrix"0"{
X \ar[r]^f \ar[dd]_p & Y \ar[d]\\
& \tilde Y' \ar[d]\\
X' \ar[r]_{f'} \ar[ur] & Y'
}
\POS-(10,-10)
\xymatrix"1"{
X_1 \ar[r]^{f_1} \ar["0"] \ar[dd]_p & Y_1 \ar[d] \ar["0"]\\
& \tilde Y_1' \ar[d] \ar["0"]\\
X_1' \ar[r]^{f_1'} \ar["0"] \ar[ur] & Y_1' \ar["0"]\\
}
\endxy
\]
\end{enumerate}
\end{lemma}

\begin{proof}
Consider the case \eqref{makesurjective}. Let $V = \coker(\pi_*X \to \pi_* Y)$. Since $Y$ is fibrant, we can find a $K$-module map $g\colon M \to Y$ for some cofibrant $K$-module $M$ realizing a section of $\pi_*Y \to V = \pi_*M$. Let $\phi\colon Y \to X'$ be a map realizing the dotted arrow, making the triangles homotopy commutative. Note that we may not be able to choose $\phi$ as a map over $B$. 

The maps $p\circ g$ and $f'\circ \phi\circ g\colon M \to Y'$ are homotopic by assumption. Choose a homotopy $H\colon \II \otimes M \to Y'$ such that $H_0 = p \circ g$ and $H_1 = f'\circ\phi\circ g$. Here $\II$ denotes the unit interval $[0,1]$. Since $f'$ is a fibration, we have a lift in the diagram
\[
\xymatrix{
M \ar[d]^{i_1} \ar[r]^{\phi\circ g} & X' \ar@{->>}[d]^{f'}\\
\II \otimes M \ar[r] \ar@{-->}[ur]^{\tilde H} & Y'.
}
\]
We now have a a map $M \to Y \times_{Y'} X'$ given by $(\tilde H_0,g)$, and it is a map over $B$ in a unique way.

Now factor the map $X \vee M \xrightarrow{(f,g)} Y$ as a trivial cofibration followed by a fibration, $\xymatrix@1@C=1pc{X \vee M \ar@{ >->}[r]^-{\sim} & \tilde X \ar@{->>}[r] & Y}$.
The map $X \vee M \to X'$ given by $(p,*)$ factors through $\tilde X$ because $X'$ is fibrant, and the commutativity of the diagram is verified.

Note that the construction of $\tilde X$ is by no means functorial. Assertion \eqref{liftsurjective} is a partial substitute for this deficiency. To show \eqref{liftsurjective}, consider the following diagram:
\[
\xymatrix{
X_1 \ar[r]^{f_1} \ar[d] & Y_1 \ar[d]\\
\tilde X \times_{X'} X_1' \ar[r] & Y \times_{Y'} Y_1'
}
\]
I claim that this diagram satisfies the assumptions of \eqref{makesurjective}. Since $\tilde X \to Y$, $X' \to Y'$, and $X_1' \to Y_1'$ are fibrations, so is the fiber product, thus the bottom map is a fibration. Moreover, $\tilde X \times_{X'} X_1'$ is fibrant as a fiber product of fibrant spaces. We need to produce a map
\[
\pi_*Y_1 \to \pi_*(\tilde X \times_{X'} X_1') = \pi_*\tilde X \times_{\pi_*X'} \pi_*X_1'
\]
making the resulting triangles commute in homotopy. A map $\pi_*Y_1 \to \pi_*X_1'$ is given by $\phi_1$; furthermore, since $\pi_*\tilde X \to \pi_*Y$ is surjective, we can find a section $\pi_*Y \to \pi_* \tilde X$ whose composite with $\pi_* \tilde X \to \pi_* X'$ is $\phi$. By the commutativity of
\[
\xymatrix{
\pi_*Y_1 \ar[r]^{\phi_1} \ar[d] & \pi_*X_1' \ar[d]\\
\pi_*Y \ar[r]^{\phi} & \pi_*X',
}
\]
we obtain a well-defined diagonal map making everything commute in homotopy.

For the assertion \eqref{makeinjective}, we apply a dual construction. Let $V = \ker(\pi_*X' \to \pi_* Y')$. Since $X'$ is cofibrant, we can find a $K$-module map $g\colon X' \to M$ for some fibrant $K$-module $M$ realizing a retraction of $\pi_*M = V \to \pi_*X'$. Let $\phi\colon Y \to X'$ be as before. 

The maps $g\circ p$ and $g \circ f\circ \phi\circ g\colon X \to M$ are homotopic by assumption. Choose a homotopy $H\colon X \to \hom(\II,M)$ such that $H_0 = g\circ p$ and $H_1 = g \circ f\circ \phi\circ g$. Since $f$ is a cofibration, we have a lift in the diagram
\[
\xymatrix{
X \ar[r] \ar@{>->}[d] & \hom(\II,M) \ar[d]^{i_1^*}\\
X' \ar[r]_-{g \circ \phi} \ar@{-->}[ur]^{H}& \hom(\{1\},M) \cong M
 }
\]
We now have a a map $X' \cup_{X} Y \to M$ given by $(g,\tilde H_0)$. Now $Y' \times M$ is a $K$-module over $B$ by means of the map $Y' \times M \to Y' \to B$. Factor $X' \xrightarrow{(f',g)} Y' \times M$ as a cofibration followed by a trivial fibration $\xymatrix@1@C=1pc{X' \ar@{ >->}[r] & \tilde Y' \ar@{->>}[r]^-{\sim} & Y' \times M}$. The map $Y \to Y' \times M$ given by $(p,*)$ factors through $\tilde Y'$ because $Y$ is cofibrant, and the commutativity of the diagram is again verified.

The proof for \eqref{liftinjective} is \emph{not} dual to \eqref{liftsurjective}. (The dual assertion would be that given $(D_1,\phi_1)$, we can find $\phi$ in $D$ compatible with $\phi_1$). It is in fact easier. Consider the diagram
\[
\xymatrix{
X_1 \ar[r] \ar[d] & Y_1 \ar[d]\\
X_1' \ar[r] & \tilde Y' \times_{Y'} Y_1'.
}
\]
This diagram satisfies the conditions of \eqref{liftinjective} because the original map $\phi$   still works. The resulting $\tilde Y_1'$ maps down to $\tilde Y$ and to $Y_1'$. Also $\pi_* \tilde Y_1' \to \pi_* Y_1'$ is injective because $\pi_* \tilde Y_1' \to \pi_*(\tilde Y' \times_{Y'} Y_1') = \pi_*\tilde Y' \times_{\pi_*Y'} \pi_*Y_1'$ is injective by construction and $\pi_*\tilde Y' \to \pi_*Y'$ is injective by assumption.
\end{proof}

\begin{prop} \label{indwe}
The ind-weak equivalences in $\Ind-\Mod_K/B$ are exactly the essentially levelwise weak equivalences. 

Similarly, the ind-pro-weak equivalences in $\Ind-\Pro-\Mod_K/B$ are the essentially levelwise weak equivalences.
\end{prop}

By the two-out-of-three property for essentially levelwise weak equivalences, we have the liberty to produce any composition of levelwise weak equivalences and ind-(pro-)isomorphisms in the proof of the proposition. 

\begin{proof}
Let $f\colon X \to Y$ be an ind-weak equivalence in $\Ind-\Mod_K/B$. We may assume $f$ is given by a levelwise map $f_s\colon X_s \to Y_s$. The condition that $f$ is an ind-weak equivalence means that for every $s$ there is an $s'>s$ and a map $\phi\colon \pi_*Y_s \to \pi_*X_{s'}$ such that the diagram
\[
\xymatrix{
\pi_* X_s \ar[r] \ar[d] & \pi_*Y_s \ar[d] \ar[dl]_{\phi}\\
\pi_*X_{s'} \ar[r] & \pi_*Y_{s'}
}
\]
commutes. By passing to a cofinal subsystem (inducing an ind-isomorphism), we may assume that there are no $t$ between $s$ and $s'$.

We will procede in two steps, first factoring $X \to Y$ as $X \to \tilde X \to Y$ where the first map is an ind-isomorphism and the second map is a levelwise map that is surjective in homotopy, and then factoring $\tilde X \to Y$ as $\tilde X \to \tilde Y \to Y$ such that $\tilde Y \to Y$ is an ind-isomorphism and $\tilde X \to \tilde Y$ is a levelwise weak equivalence.

First, using functorial cofibrant-fibrant replacement, we may assume that all $X_s$ and $Y_s$ are cofibrant and fibrant in $\Mod_K/B$. By another functorial factorization, we may assume that $X_s \to Y_s$ is a fibration in $\Mod_K$ for all $s$. All these operations induce levelwise weak equivalences.

Applying Lemma~\ref{levelwisewelemma}\eqref{makesurjective} inductively, we obtain a diagram
\[
\xymatrix{
X_s \ar[r] \ar[d] &\tilde X_s \ar[d] \ar[dl] \ar[r]^{\tilde f_s} & Y_s \ar[d]\\
X_{s'} \ar[r] \ar[d] & \tilde X_{s'} \ar[d] \ar[dl] \ar[r]^{\tilde f_{s'}} & Y_{s'} \ar[d]\\
X_{s''} \ar[r] \ar[d] & \tilde X_{s''} \ar[d] \ar[dl] \ar[r]^{\tilde f_{s''}} & Y_{s''}\ar[d]\\
\vdots & \vdots & \vdots
}
\]
with maps $\tilde f$ which are surjective in $\pi_*$. Thus, we have found an ind-isomorphism $X \to \tilde X$ and a levelwise surjective map $\tilde X \to Y$.

Now let us assume that $X \to Y$ is a levelwise cofibration (in $\Mod_K$) of fibrant-cofibrant objects in $\Mod_K/B$, which is levelwise surjective in $\pi_*$. Then, arguing as before but using Lemma~\ref{levelwisewelemma}\eqref{makeinjective}, we get a factorization of $X \to Y$ as a map $X \to \tilde Y$ which is a levelwise isomorphism in $\pi_*$, followed by an ind-isomorphism $\tilde Y \to Y$.

The proof for $\Ind-\Pro-\Mod_K/B$ is very similar. Without loss of generality by Lemma~\ref{levelwise}, let $f\colon X \to Y$ be a levelwise map, where $X,\; Y\colon Q_{S,\alpha} \to \Mod_K/B$ are functors. We assume the we have prepared $f$ by levelwise cofibrant/fibrant replacement as before, so that Lemma~\ref{levelwisewelemma} is applicable when we need it.

Since $f$ is assumed to be an ind-pro-weak equivalence, this means that for every $s$ there is an $s'>s$ such that for every $t'\in \alpha(s')$ there is a $t \in \alpha(s)$, $t < \alpha(s < s')(t')$, and a map $\phi\colon \pi_*Y_s^t \to \pi_*X_{s'}^{t'}$ such that the diagram
\[
\xymatrix{
\pi_* X_s^t \ar[r] \ar[d] & \pi_*Y_s^t \ar[d] \ar[dl]_{\phi}\\
\pi_*X_{s'}^{t'} \ar[r] & \pi_*Y_{s'}^{t'}
}
\]
commutes. As before, by passing to a cofinal subsystem, we may assume that $s'$ is a direct successor of $s$. In the first step, we apply Lemma~\ref{levelwisewelemma}\eqref{makesurjective} and \eqref{liftsurjective} to produce a factorization $X_s^t \to \tilde X_s^t \to X_{s'}^{t'}$ with $X_s^{t'} \to Y_s^{t}$ surjective in homotopy. We cannot simply do this for every $s$ and $t$ because the construction of $\tilde X$ in Lemma~\ref{levelwisewelemma} is not functorial. Fix $s$ and assume $\tilde X_s^t$ has been constructed for some $t'<t_1'$. Consider the diagram
\[
\xy
\xymatrix"0"{
X_s^t \ar[r]^f \ar[d] & Y_s^t \ar[dd]\\
\tilde X_s^t \ar[ur] \ar[d]\\
X_{s'}^{t'} \ar[r]_{f'} & Y_{s'}^{t'}
}
\POS-(10,-10)
\xymatrix"1"{
X_s^{t_1} \ar[r]^{f_1} \ar["0"] \ar@{-->}[d] & Y_s^t \ar[dd] \ar["0"]\\
\tilde X_s^{t_1} \ar@{-->}["0"] \ar@{-->}[ur] \ar@{-->}[d]\\
X_{s'}^{t_1'} \ar[r]^{f_1'} \ar["0"] & Y_{s'}^{t_1'} \ar["0"]\\
}
\endxy
\]
By Lemma~\ref{levelwisewelemma}\eqref{liftsurjective}, we can find $\tilde X_s^{t_1}$ and maps as indicated in the diagram. Proceeding inductively, we obtain a pro-object $\tilde X_s$ for all $s \in S$. For varying $s$, these assemble to an ind-pro-object by means of the maps
\[
\tilde X_s^{t} \to X_{s'}^{t'} \to \tilde X_{s'}^{t'}.
\]
Furthermore, the commutative diagram
\[
\xymatrix{
X_s^{t_1} \ar[r] \ar[d] & \tilde X_s^{t_1} \ar[d] \ar[dl]\\
X_{s'}^{t'} \ar[r] & \tilde X_{s'}^{t'}
}
\]
shows that $X \to \tilde X$ is an ind-pro-isomorphism. By construction, $\tilde X_s^t \to Y_s^t$ is surjective in homotopy.

We leave the dual construction of $\tilde Y \to Y$ and of a levelwise weak equivalence $\tilde X \to \tilde Y$ to the reader.
\end{proof}

\section{Independence of the total space} \label{indeptotal}

The aim of this section is to prove Theorem~\ref{universalconvergence}.

Fix a multiplicative homology theory $K$. For the results of this section, $K$ does not need to be a field (i.e. have K\"unneth isomorphisms for any two spaces). We define the categories
\[
\T = \Ind-\Pro-\Mod_K \text{ and } \A = \Ind-\Pro-\Mod_{K_*}.
\]
For a CW-complex $X$, let $\Fin(X)$ be the directed set of finite subcomplexes of $X$. Denote by $K^{\fib} \colon (\Top\times\Top)/B \to \T$ the functor with
\[
K^{\fib}(X_1 \to B \leftarrow X_2) = \{K[F']\}_{F' \in \Fin(\holim(X_1 \to B \leftarrow X_2))}
\]
where $K[F'] \in \Pro-\Mod_K$ as an object indexed over the one-point category, or, according to taste, as a constant tower.

Similarly, define a functor $\KK \colon (\Top\times\Top)/B \to \T$ by
\[
\KK(X_1 \to B \leftarrow X_2) = \{\Tot^s K[C_B(X_1',X_2')]\}_{s \geq 0, X_1'\in \Fin(X_1),\; X_2' \in \Fin(X_2)}.
\]

We write $K^{\fib}_* = \pi_*$ and $\KK_* = \pi_*\KK\colon K^{\fib}\colon (\Top\times\Top)/B \to \A$.

The map $\Phi\colon K[\holim X_1 \to B \leftarrow X_2] \to \{\Tot^s K[C_B(X_1,X_2)]\}_s$ extends to a natural transformation $\Phi\colon K^{\fib} \to \KK$ as follows: if $F'$ is a finite subcomplex of $\holim(X_1 \to B \leftarrow X_2)$ then its images in $X_1$ and $X_2$ are finite subcomplexes, and we get a comparison map
\[
K[F'] \to \{\Tot^s K[C_B(\im(F' \to X_1),\im(F' \to X_2))]\}_{s \geq 0}
\]
These assemble to a map in $\T$. 

\begin{defn}We call a diagram $X_1 \to B \leftarrow X_2$ of spaces \emph{ind-pro-constantly convergent}, or, more briefly, \emph{Ind-convergent}, if $\pi_*\Phi(X_1 \to B \leftarrow X_2)$ is an isomorphism in $\A$. 

We call a map $X \to B$ Ind-convergent if $X \to B \leftarrow Y$ is Ind-convergent for every $Y \to B$, and we call a space $B$ Ind-convergent if every homotopy pullback diagram $X_1 \to B \leftarrow X_2$ is Ind-convergent.
\end{defn}

\begin{remark}
By saying that $\pi_* \Phi$ is an isomorphism, we really mean that $\Phi$ induces an isomorphism after applying the graded group valued functor $\prod_i \pi_i$, which a priori is stronger than requiring that for every $k$, $\pi_k$ induces an isomorphism. However, since in the context of this paper (if not this section) all homotopy groups are homotopy groups of $K$-modules for periodic theories $K$, the two notions coincide.
\end{remark}

Note that if $X_1$ and $X_2$ are finite CW-complexes, $\KK(X_1 \to B \leftarrow X_2)$ is ind-constant, but $K^{\fib}(X_1 \to B \leftarrow X_2)$ is ind-constant if and only if $K_*(F)$ is finite. Thus Ind-convergence of $B$ implies in particular that the fiber of a fibration with total space a finite CW-complex has finite $K$-homology. In this situation, Ind-convergence is the same as pro-constant convergence.

\begin{example}
This example shows that Ind-convergence is weaker than pro-con\-stant convergence. Let $B=\S^1$, $E_i = \S^1$, and $p_i\colon E_i \to B$ multiplication by $i$. Thus, $F_i$ is the discrete space with $i$ points. The $H\Z$-based Eilenberg-Moore spectral sequence for $E_i \to B$ is $i$-convergent; more specifically,
\[
E^2_{**} = \Z[x] \otimes \bigwedge(y),
\]
where $x$ is in bidegree $(-1,1)$ and $y$ is in bidegree $(0,1)$, and we have differentials
\[
d^i(y) = x^i.
\]
This shows that for $E = \coprod_i E_i $, $E \to B$ cannot be pro-constantly convergent because there are differentials of arbitrary length. However, $E \to B$ is Ind-convergent because $\coprod_{i \leq n} E_i \to B$, which is $n$-convergent, constitutes a cofinal subsystem of the finite sub-CW-complexes of $E$.
\end{example}

\begin{thm} \label{universalconvergencehelper}
Let $Y \to B \leftarrow *$ be an Ind-convergent map for some $Y \to B$, where $K_*$ is a graded field. Then $Y \to B$ is Ind-convergent.
\end{thm}

A cohomological version of this theorem (in terms of pro-constant convergence) was proven in \cite{hodgkin:emss,jeanneret-osse:emss} under some cohomological finiteness conditions on $K^*(X)$. More strongly, \cite{seymour:emss} claims that the cohomological finiteness condition is not necessary if $B$ is a finite-dimensional CW-complex, but the proof seems to contain mistakes. Our formulation does not require any such restriction; however we need the rather strong assumption of Ind-convergence to begin with.

Theorem~\ref{universalconvergence} is an immediate corollary.

Fix a map $Y \to B$, and abbreviate $\KK(X \to B \leftarrow Y)$ as $\KK(X)$ and $K^{\fib}(X \to B \leftarrow Y)$ as $K^{\fib}(X)$.

\begin{lemma}\label{homologytheories}
$K^{\fib}_*(X)$ and $\KK_*(X)$ are homology theories on the category $(\Top/B)_*$ of sectioned spaces over $B$ with values in $\A$ in the sense of Dold \cite{dold:chern-classes}. This means: in addition to the usual axioms for a homology theory $h$ on $(\Top/B)_*$ (long exact sequence, excision), the following two axioms are satisfied:
\begin{itemize}
\item[(CYL)] For any $X \to B$, $p\colon X \times [0,1] \to B$, $h_*(X \times \{0\}) \to h_*(X \times [0,1])$ is an isomorphism;
\item[(EXC)] If 
\[
\xymatrix{ A \ar[r] \ar[d] & X_1 \ar[d] \\ X_2 \ar[r] & X}
\]
is a pushout square in $\Top$ and $p\colon X \to B$ is a map then we have an exact Mayer-Vietoris sequence.
\end{itemize}
\end{lemma}

\begin{remark}
The additional axioms have an analog for $G$-spaces. In general, if a $G$-equivariant map $X \to Y$ is also a homotopy equivalence, it need not induce an isomorphism on some $G$-equivariant homology theory because the homotopy inverse need not be $G$-equivariant. It does so, however, if the homology theory only depends on the \emph{homotopy} fixed point or orbits for some subgroup; for example, Borel homology satisfies the analogue of (CYL) and (EXC). 
\end{remark}

\begin{proof}[Proof of Lemma~\ref{homologytheories}]
Let us first consider $K^{\fib}$. The functor
\[
X \mapsto K_*(F) = \colim K^{\fib}_*(X)
\]
is a homology theory by \cite[3.4]{dold:chern-classes}; furthermore a sequence in $\Ind-\Mod_{K_*}$ is exact if and only if its colimit sequence is exact. This shows that we have long exact sequences. To verify the wedge axiom, we show that $K^{\fib}$ maps filtered hocolimits to colimits. Note that
\[
\holim(\hocolim_i X_i \to B \leftarrow Y) \simeq \hocolim_i \holim(X_i \to B \leftarrow Y),
\]
so that we only need to see that
\[
\{K_*F\}_{F' \in \Fin(\hocolim_i F_i)} \cong \colim_i  \{K_*(F_i')\}_{F_i' \in \Fin(F_i)} \cong \{K_*(F_i')\}_{i, F_i' \in \Fin(F_i)}
\]
where the colimit is taken in $\Ind-\Mod_{K_*}$ and the last isomorphism is its definition in any ind-category. Since any finite sub-CW-complex of a hocolim is already contained in an $F_i$, the two indexing systems are mutually cofinal, and the isomorphism is shown. Axiom (CYL) is clearly satisfied because the homotopy fibers of $X$ and of $X \times [0,1]$ are homotopy equivalent; finally, (EXC) is satisfied because the homotopy pullback functor sends pushout squares in a total space to pushout squares.

Now consider $\KK$. The functor $X \mapsto \pi_*(\Tot^s K[C_B(X,Y)])$ is a homology theory by construction (as an iteration of taking cofibers and smash products with B and Y of the suspension spectrum of X), and levelwise exact sequences induce exact sequences in $\Pro-\Mod_{K_*}$, thus $X \mapsto \{\pi_*(\Tot^s K[C_B(X,Y)])\}_s$ is a homology theory with values in $\Pro-\Mod_{K_*}$ (not satisfying the wedge axiom!) and thus induces a homology theory on \emph{finite} CW-complexes over $B$. 

Now let $U \to V$ be a map in $(\Top/B)_*$, we may assume an inclusion of CW-complexes. For any $V' \in \Fin(V)$, let $U' = U \cap V' \in \Fin(U)$. Then the sequence $\KK_*(U) \to \KK_*(V) \to \KK_*(V/_B U)$ has the level representation
\[
\{\KK_*(U') \to \KK_*(V') \to \KK_*(V'/_B U')\}_{V' \in \Fin(V)}.
\]
Since this sequence is levelwise exact by the above, it is exact in $\A$. The wedge axiom is also satisfied:
\[
\KK(\hocolim X_i) = \{\KK(X')\}_{X' \in \Fin(\hocolim(X_i)} = \{\KK(X')\}_{i, \; X' \in \Fin(X_i)} = \colim \KK(X_i).
\]
The axioms (CYL) and (EXC) are satisfied levelwise in $\Pro-\Mod_K$ whenever the total spaces are finite CW-complexes, and an argument similar to the one above for exact sequences shows that they hold essentially levelwise in $\A$.  
\end{proof}

The added value of the axioms (CYL) and (EXC) is that a natural transformation of homology theories satisfying the axioms is a natural isomorphism if and only if it is an isomorphism on points \cite[Theorem~4.1]{dold:chern-classes}:

\begin{thm}\label{doldthm}
Let $\Phi\colon h \to h'\colon (\Top/B)_* \to A$ be a natural transformation of homology theories, where $A$ is some abelian category. Then $\Phi(X \to B)$ is an isomorphism for all $X \to B$ iff $\Phi(* \to B)$ is an isomorphism for all points in $B$. \qed
\end{thm}

\begin{proof}[Proof of Theorem~\ref{universalconvergencehelper}]
Putting $h=K^{\fib}$, $h'=\KK$ in Theorem~\ref{doldthm}, using Lemma~\ref{homologytheories}, we obtain that $\Phi(X \to B)$ is an isomorphism.
\end{proof}

\section{Transitivity of convergence} \label{transitivitysection}

The principal aim of this section is to prove Theorem~\ref{fibrationconvergence}. This will follow rather easily from the following, more general result.

\begin{thm} \label{transitivity}
Let $F_1 \to X \xrightarrow{\pi_1} B_1$, $F_2 \to X \xrightarrow{\pi_2} B_2$ be two fibrations of connected spaces. Denote by $F$ the fiber of $F_1 \to B_2$, which is the same as the fiber of $F_2 \to B_1$. Let $K$ be a homology theory which is a field, and assume that $X \to B_i$ are Ind-convergent for $i=1,\;2$. Then the fibration $F_1 \to B_2$ is Ind-convergent if and only if the fibration $F_2 \to B_1$ is Ind-convergent.
\end{thm}

The following corollary is the special case of $B=B_1$, $B_2 = X$ and explains why I call this result a transitivity property.

\begin{corollary}\label{transitivecor}
Let $F \to X \to B$ be a fibration sequence such that $X \to B$ and $F \to X$ are Ind-convergent. Then $B$ is Ind-convergent. \qed
\end{corollary}

We are now in a position to derive Theorem~\ref{fibrationconvergence} from Theorem~\ref{transitivity}.

\begin{corollary}[Theorem~\ref{fibrationconvergence}]
Let $F \to Y \to X$ be a fibration with $F$ and $X$ Ind-convergent. Then so is $Y$.
\end{corollary}
\begin{proof}
Consider the diagram
\[
\xymatrix{
\Omega X \ar[r] \ar[d] & {*} \ar[r] \ar[d] & X \ar@{=}[d]\\
F \ar[r] \ar[d] & Y \ar[r] \ar@{=}[d] & X \ar[d]\\
Y \ar@{=}[r] & Y \ar[r] & {*}
}
\]
Since $X$ is Ind-convergent, both upper rows are Ind-convergent, and the middle vertical row is $0$-convergent because the fibration is trivial. By Theorem~\ref{transitivity}, $F \to Y$ is Ind-convergent (to $K_*(\Omega X)$). Applying Corollary~\ref{transitivecor} to $\Omega X \to F \to Y$, we find that $Y$ is Ind-convergent.
\end{proof}

The technical heart of Theorem~\ref{transitivity}, and in fact the raison d'\^etre for all of Section~\ref{indprosection}, is the following lemma.

\begin{lemma}\label{cobarindequivalence}
Let $f\colon X \to Y$ be a map in $\Ind-\Pro-\Mod_K/B$ such that $\pi_*f\colon \pi_*X \to \pi_*Y$ is an isomorphism in $\A$. Then $\Tot^\bullet C_B(X) \to \Tot^\bullet C_B(Y)$ also induces an isomorphism in $\A$.
\end{lemma}
\begin{remark}
By the Five Lemma, it is obvious that $f$ induces an isomorphism in $\Pro-\Ind-\Pro-\Mod_{K_*}$, but it is not obvious that it lifts to $\Ind-\Pro-\Mod_{K_*}$.
\end{remark}
\begin{proof}
By assumption, $f\colon X \to Y$ is an ind-pro-weak equivalence in the category $\Ind-\Pro-\Mod_K/B$. By Proposition~\ref{indwe}, we may assume that there is a directed set $S$, a functor $\alpha\colon S^{\op} \to \{\text{inverse sets}\}$, functors $M ,\; N\colon \QQ_{S,\alpha} \to \Mod_K/B$, and a commutative diagram
\[
\xymatrix{
X \ar[r]^{K \sm f} \ar[d]_{\Ind-\Pro-\cong}^{\alpha_X} & Y \ar[d]^{\Ind-\Pro-\cong}_{\alpha_Y}\\
M \ar[r]^{\tilde f} & N
}
\]
where the vertical maps $\alpha_X$ and $\alpha_Y$ are ind-pro-isomorphisms and $\tilde f$ is a levelwise weak equivalence. Now note that the functor
\begin{align*}
\Tot^\bullet C_B\colon \Mod_K/B &\to \Pro-\Mod_K\\
X &\mapsto \{\Tot^s C_B(X)\}_{s \geq 0}
\end{align*}
extends to a functor
\begin{align*}
\Ind-\Pro-\Mod_K/B &\xrightarrow{\Ind-\Pro-\Tot^\bullet C_B} \Ind-\Pro-\Pro-\Mod_K\\
&\xrightarrow{\Ind-D} \Ind-\Pro-\Mod_K,
\end{align*}
where $D\colon \Pro-\Pro-\C \to \Pro-\C$ is the tautological functor of Remark~\ref{proproremark}.

We thus obtain isomorphisms
\[
\Tot^\bullet C_B(X) \cong \Tot^\bullet C_B(M) \quad \text{and} \quad \Tot^\bullet C_B(Y) \cong \Tot^\bullet C_B(N)
\]
as well as a levelwise weak equivalence $C_B(\tilde f)\colon C_B(M) \to C_B(N)$. This induces a levelwise weak equivalence of total towers in every degree and thus an ind-pro-isomorphism
\begin{multline*}
\pi_*\Tot^\bullet C_B(X) \cong_{\Ind-\Pro} \pi_*\Tot^\bullet C_B(M) \\
\cong \Tot^\bullet C_B(N) \cong_{\Ind-\Pro} \pi_*\Tot^\bullet C_B(Y).
\end{multline*}
\end{proof}

Consider the diagram of fibrations
\begin{equation}\label{transitivitydiagram}
\xymatrix{
F \ar[r] \ar[d] & F_1 \ar[d] \ar[r] & B_2 \ar@{=}[d]\\
F_2 \ar[r] \ar[d] & X \ar[r]^{\pi_2} \ar[d]_{\pi_1} & B_2 \ar[d]^{\top}\\
B_1 \ar@{=}[r] & B_1 \ar[r]^{\top} & {*}
}
\end{equation}
To compare the two Eilenberg-Moore spectral sequences abutting to the $K$-homo\-logy of $F$, we construct a bicosimplicial space
\[
C^{st}(X) = X \times B_1^s \times B_2^t
\]
We have
\[
\Tot C^{\bullet t}(X) = F_1 \times B_2^t \quad \text{and thus} \quad \Tot \Tot C^{\bullet\bullet}(X) = F;
\]
on the other hand we have
\[
\Tot C^{s \bullet}(X) = B_1^s \times F_2.
\]

We denote the ``horizontal'' total space $\Tot \{C^{st}\}_s$ by $\Tot_h C^{\bullet\bullet}$ and the ``vertical'' total space by $\Tot_v C^{\bullet\bullet}$.

\begin{lemma}\label{innerindiso}
Given a diagram as in \eqref{transitivitydiagram} with $X \to B_1$ Ind-convergent, there is an isomorphism
\[
\KK(F_1 \to B_2) \xrightarrow{\cong} \KK(X \to B_1 \times B_2).
\]
\end{lemma}
\begin{proof}
Since $X \to B_1$ is Ind-convergent, the map
\[
\{K[F]\}_{F \in \Fin(F_1)} = K^{\fib}(X \to B_1) \to \KK(X \to B_1) = \{ \Tot^\bullet K[C_{B_1}(X')]\}_{X' \in \Fin(X)}.
\]
is an ind-pro-weak equivalence as well as a map over $B_2$. By Lemma~\ref{cobarindequivalence}, it induces an ind-pro-weak equivalence
\begin{align*}
\KK(F_1 \to B_2) = &\{\Tot^\bullet K[C_{B_2}(F')]\}_{F' \in \Fin(F_1)}\\
 \to & \{\Tot^\bullet C_{B_2}(\Tot^\bullet C_{B_1}(X'))\}_{X' \in \Fin(X)}\\
\cong &\{\Tot^\bullet \Tot_h^\bullet K[C^{\bullet\bullet}(X')]\}_{X' \in \Fin(X)}.
\end{align*}
Since the diagonal $\N \hookrightarrow \N \times \N$ is cofinal, the right hand side is
\[
\{\Tot^\bullet \Tot_h^\bullet K[C^{\bullet\bullet}(X')]\} \cong \{\diag \Tot^\bullet \Tot_h^\bullet K[C^{\bullet\bullet}(X')]\} \cong \KK(X \to B_1 \times B_2).
\]
\end{proof}

\begin{proof}[Proof of Theorem~\ref{transitivity}]
From Diagram~\ref{transitivitydiagram}, we obtain a diagram in $\T$ from the various comparison maps
\[
\xymatrix@C=2cm{
K^{\fib}(F_2 \to B_1) \ar@{=}[d] \ar[rr]^-{\text{convergence of $F_2 \to B_1$}} && \KK(F_2 \to B_1) \ar[d]^{\text{Lemma~\ref{innerindiso} for $X \to B_2$}}_{\sim}\\
K^{\fib}(X \to B_1 \times B_2) && \KK(X \to B_1 \times B_2)\\
K^{\fib}(F_1 \to B_2) \ar@{=}[u] \ar[rr]^-{\text{convergence of $F_1 \to B_2$}} && \KK(F_1 \to B_2). \ar[u]_{\text{Lemma~\ref{innerindiso} for $X \to B_1$}}^{\sim}\\
}
\]
The two-out-of-three property in this diagram finishes the proof of the theorem.
\end{proof}

\section{The Hopf coring for Morava $K$-theory of Eilenberg-Mac Lane spaces} \label{cohopf}

In this section, we completely analyze the $K(n)$-based Eilenberg-Moore spectral sequence for path-loop fibrations on mod-$p$ Eilenberg-Mac~Lane spaces, where $p$ is the same prime as the characteristic of $K(n)_*$. As $n \geq 1$ is fixed in this section, we will abbreviate $K(n)$ by $K$. Throughout this section, $p$ is assumed to be odd (an assumption made in the crucial input \cite{ravenel-wilson:morava}, and also necessary to ascertain that $K(n)$ is a homotopy commutative ring spectrum). It is likely that the convergence result also holds for $p=2$.

In order to understand the $K$-based Eilenberg-Moore spectral sequence for the Eilenberg-Mac Lane spaces $\EMH_i = K(\Z/p,i)$, it will be necessary to understand the algebra structure of $K^*(\EMH_*)$ quite well. Ravenel and Wilson study the structure of $K_*(\EMH_*)$ in \cite{ravenel-wilson:morava}, and all the necessary information can be extracted from that paper. However, they use the \emph{bar spectral sequence}, which is known to converge, to compute $K_*\EMH_{i+1}$ from $K_*\EMH_i$. In our case, we are interested in the Eilenberg-Moore spectral sequence, and from a knowledge of $K_*\EMH_{i+1}$ and $K_*\EMH_i$ we want to conclude that the spectral sequence converges for $i<n$. Since we know what the answer \emph{should} be, we could guess what the differentials would have to be, and that guess is in fact correct, but unfortunately provides no proof for convergence. Thus we need to study $E_*\EMH_i$ for other homology theories $E \neq K$, for which we know the EMSS converges, and then compare it to the $K$-based EMSS, deriving the differentials there in a rigid way. By the various multiplicative properties of the EMSS as described below, we can restrict ourselves to computing the $E$-based EMSS for $\EMH_1$ and $E=k(n)$, connective Morava $K$-theory, to derive the differentials in all other cases for $K$.

Since $K$ is a graded field, we have a K\"unneth isomorphism
\[
K_*(\EMH_r \times \EMH_s) \cong K_*(\EMH_r) \otimes_{K_*} K_*(\EMH_s)
\]
for all $r$ and $s$, and thus $K_*(\EMH_r)$ is a coalgebra, and therefore the gadget $K_*(\EMH_*)$ obtains the structure of a \emph{Hopf ring}, i.~e., it is a ring object in the category of coalgebras. For a survey on this kind of algebraic structure, consult \cite{wilson:hopf-rings}. The important data here are operations $\Psi$, $+$, $*$, $\circ$, where
\begin{align*}
\Psi\colon  K_s(\EMH_r) &\to \left(K_*(\EMH_r) \otimes_{K_*} K_*(\EMH_r)\right)_s;\\
+\colon K_s(\EMH_r) \otimes K_s(\EMH_r) &\to K_s(\EMH_r);\\
*\colon  K_{s'}(\EMH_r) \otimes K_{s''}(\EMH_r) &\to K_{s'+s''}(\EMH_r);\\
\circ:  K_{s'}(\EMH_{r'}) \otimes K_{s''}(\EMH_{r''}) &\to K_{s'+s''}(\EMH_{r'+r''}).
\end{align*}
 
The coproduct and addition are the usual maps in homology; the  $*$-product is the ``additive'' product coming from the (infinite) loop space structure of $\EMH_r$; and the $\circ$-product is the ``multiplicative'' product coming from the ring spectrum map $\EMH_r \times \EMH_{r'} \to \EMH_{r+r'}$. As usual, $*$ distributes over $+$, 
but there is a second layer of distributivity; namely, in the Sweedler notation $\Psi(a) = \sum_{(a)}a' \otimes a''$,
\[
a \circ (b*c) = \sum_{(a)} (a' \circ b) * (a'' \circ c).
\]
By convention, we give $\circ$ operator precedence over $*$, so that we could write the summand in the above formula without parentheses.

Both products are $\Psi$-comodule maps. There are, of course, a number of other structural maps corresponding to units, counits, and coinverses. We denote by $[1]: \pi_0(\EMH_0) \to K_0(\EMH_{0})$ the image of the unit under the Hurewicz homomorphism, which is the unit for the $\circ$-product. Similarly, denote by  $[0]_r \in K_0(\EMH_r)$ are the units for the $*$-products in degree $r$.

\begin{notation}
In our computations, we will need to deal with algebras, coalgebras, and Hopf algebras over $\F_p$ or $K_*$. We adopt the following standard notation:
\begin{description}
\item[$P(x)$] is the Hopf algebra whose underlying algebra is the polynomial algebra on $x$ and whose underlying coalgebra is the divided polynomial coalgebra. We denote the standard additive generators by $x^i$, as usual.
\item[$P_k(x)$] is the quotient of $P(x)$ whose underlying algebra is the truncated polynomial algebra $P(x)/(x^{p^k})$.
\item[$\Gamma(x)$] is the Hopf algebra dual of $P(x)$: its underlying algebra is the divided polynomial algebra, and its underlying coalgebra is the tensor coalgebra. We denote the standard additive generators by $x_i$.
\item[$\Gamma_k(x)$] is the sub-Hopf algebra of $\Gamma(x)$ on the generators $x_i$ ($0 \leq i < p^k$).
\item[$\exterior(x)$] is the exterior Hopf algebra on a primitive generator.
\item[$R_u(x)$] is the Hopf algebra $P(x)/(x^p-ux)$ (where $u$ is a unit) with $x$ primitive.
\end{description}
We will use the convention that $x^{(i)} = x^{*p^i}$ and $x_{(i)} = x_{p^i}$.
If $x \in K_s(\EMH_r)$, we will write $|x|=(s,r)$.
\end{notation}

The following lemma is basic multiplicative homological algebra:
\begin{lemma}\label{basictor}
If $A$ is a (graded) commutative algebra over a field $k$ of characteristic $p$, $\Tor_{**}(A) =_{\text{def}} \Tor_{**}^A(k,k)$ is a commutative and cocommutative Hopf algebra. In particular,
\begin{itemize}
\item $\Tor_{**}(\exterior(y)) \cong \Gamma(\sigma y);$
\item $\Tor_{**}(P(x)) \cong \exterior(\sigma x);$
\item $\Tor_{**}(P_n(x)) \cong \exterior(\sigma x) \otimes \Gamma(\phi x).$
\end{itemize}
Here and in the following, $\sigma x$ denotes the suspension, i.~e. the element in $\Tor_1$ which is represented in the bar resolution by $[a]$, and $\phi x$ denotes the ``transpotence'' element in $\Tor_2$, which is represented in the bar resolution by any of the classes $[ x^i \mid x^j ]$ with $i,\;j \geq 1$, $i+j=p^n$ (up to units).

Dually, if $C$ is a cocommutative coalgebra over $k$, $\Cotor_{**}(C) =_{\text{def}} \Cotor_{**}^C(k,k)$ is a commutative and cocommutative Hopf algebra as well, with
\begin{itemize}
\item $\Cotor_{**}(\exterior(y)) \cong P(\sigma y);$
\item $\Cotor_{**}(\Gamma(x)) \cong \exterior(\sigma x);$
\item $\Cotor_{**}(\Gamma_n(x)) \cong \exterior(\sigma x) \otimes P(\phi x).$
\end{itemize}
Again, $\sigma x$ denotes the dual of the suspension and $\phi x$ denotes the ``cotranspotence'' element in $\Cotor_2$, which is represented in the bar resolution by $\sum_{i+j=n\atop i,\;j \geq 1} u_i [ x^i \mid x^j ]$ for some units $u_i$ we do not care about.
\end{lemma}

\subsection{Getting started: $K_*\EMH_1 \to K_*\EMH_0$}

We start our computation by studying the differentials of the EMSS for the path-loop fibration on $\EMH_1 = K(\Z/p,1)$ in connective Morava $K$-theory $k$. We have, as coalgebras,
\begin{align}
E_*(\EMH_0) &= E_*[\Z/p] & \text{for any $E$}\\
k^{\ev}_*(\EMH_1) &= \Gamma_n(a) \label{andefinition}\\
k^{\odd}_*(\EMH_1) &= \langle y_1,y_2,\dots\rangle/(v_n) & \text{with $|y_i|=2i-1$}
\end{align}
The result for $k$ is easily computed with the Atiyah-Hirzebruch spectral sequence. Since $k$ is not a graded field, $k_*(\EMH_1)$ need not be a coalgebra, and indeed is not. But we can compute the cobar spectral sequence modulo the Serre class of $v_n$-torsion groups, so that we have a K\"unneth isomorphism again and an isomorphism of Hopf algebras
\[
E^2_{s,t} = \Cotor_{s,t}(\Gamma_n(x)) = \exterior(\sigma a) \otimes P(\phi a) \quad \text{(Lemma~\ref{basictor})}
\]

The only way this can converge is to have a differential
\[
d^{2p-1}(v_n^p \sigma a) = (\phi a)^{(1)}.
\]
Inverting $v_n$, we derive the same differential in the $K$-based EMSS, and we thus have convergence there, too, with
\[
E^\infty_{s,t} = E^{2p}_{s,t} = P_1(\phi a),
\]
where $v_n^{-1}\phi a$ represents $[1]-[0] \in K_0(\EMH_0)$, and we have a $*$-multiplicative extension $(v_n^{-1}\phi a)^{(1)} \dashrightarrow v_n^{-1}\phi a$.  In fact, we have $2p$-convergence since there are no longer differentials and the filtration in $E^\infty$ is bounded by $2p$.

\subsection{Morava $K$-theory of Eilenberg-Mac~Lane spaces a.k.a. Automorphisms of Siegel domains}

Unfortunately, we will need to juggle around with multi-indices quite a bit. A multi-index $I$ is an $n$-tuple $(i_0,i_1,\dots,i_{n-1})$ with $i_\nu \in \{0,1\}$.
For such multi-indices we will need some operators, functions, and constructions, which we assemble in the following definition. The reader is advised to skip it and refer back to it when the notation is used. Beware that our $sI$ is what in \cite{ravenel-wilson:morava} and \cite{wilson:hopf-ring-morava} would be called $s^{-1}I$, and our usage of multi-indices differs from \cite{ravenel-wilson:morava} but agrees with \cite{wilson:hopf-ring-morava}.

\begin{defn}$\;$

\begin{itemize}
\item \textbf{Constructions.} Denote by $\Delta_k$ the index 
$(\delta_{lk})_l$, where $0 \leq l \leq n-1$.

Denote by $\Delta[k]$ the index $(1,\dots,1,0,\dots,0)$ ($k$ copies 
of ones) and
$\nabla[k]$ for the index $(0,\dots,0,1,\dots,1)$ ($k$ copies 
of ones).
\item \textbf{Operations.} For a multi-index $I=(i_0,i_1,\dots,i_{n-1})$, let 
$s(I)$ denote the shift $(i_1,i_2,\dots,i_{n-1},0)$.

Also denote by $cI$ the cyclic permutation $(i_1,i_2,\dots,i_{n-1},i_0)$.
\item \textbf{Functions.} Let $t_\epsilon(I)+1$ (the number of 
trailing $\epsilon$) denote the smallest $k$ with $i_{n-k}\neq\epsilon$, 
or $\infty$ if $I=(1-\epsilon,\dots,1-\epsilon)$. Denote by 
$l_\epsilon(I)$ (the number of leading $\epsilon$) the smallest $k$ 
such that $i_k=0$.
\end{itemize}
\end{defn}

We write
\[
a^I = a_{(0)}^{\circ i_0} \circ a_{(1)}^{\circ i_1} \circ a_{(n-1)}^{\circ i_{n-1}}
\]
where the $a_{(i)} \in K_*(\EMH_1)$ are defined as in \eqref{andefinition}. Thus, $|a^I| = (\sum_{\nu} i_\nu p^\nu,\sum I)$.

\begin{thm}[Ravenel-Wilson] \label{wilsonthm}
In terms of the classes defined above, we have an isomorphism of $K_*$-algebras
\[
K_*(\EMH_*) \cong \bigotimes_{I \atop i_0 = 0} P_{t_1(I)+1}(a^I) \otimes R_u(a^{\Delta[n]}),
\]
for some unit $u \in K_*^\times$, where $I=(i_0,\dots,i_{n-1})$ runs through all multi-indices.

The coproduct is completely determined by stating that
\begin{equation}\label{standardcoprod}
\Psi(a_i) = \sum_{j=0}^i a_{i-j} \otimes a_j.
\end{equation}
\end{thm}

The classes $a^I$ for $i_0=1$ do not appear as generators in Theorem~\ref{wilsonthm}, but they are nonzero and thus can be expressed in terms of the generators. This computation is a reformulation of Ravenel-Wilson's.

\begin{lemma}\label{normalizedrep}
Let $m=l_1(I)$, $I \neq \Delta[n]$. Then
\[
a^I = (-1)^{m\sum I} \left(a^{c^mI}\right)^{(m)}
\]
\end{lemma}

\begin{corollary}
As $K_*$-modules,
\[
K_*(\EMH_*) = \smashop\bigotimes_I P_1(a^I) = M_{**}.
\]
The coalgebra structure is given by \eqref{standardcoprod} and the fact that $\Psi$ is an algebra  morphism with respect to $\circ$ and also with respect to the above algebra structure for the multiplication $*$. \label{coprodlinear}

In particular, if
\[
\Psi^p: K_*(\EMH_r) \to \left(K_*(\EMH_r)\right)^{\otimes p}
\]
denotes $p$-fold comultiplication, then
\[
\Psi^p(a^I) = \begin{cases}
0; & \text{if $i_0 =1$}\\
\left( a^{sI}\right)^{\otimes p} + \text{decomposables}; & \text{otherwise.}
\end{cases}
\]
\end{corollary}
\begin{proof}

Using Lemma~\ref{normalizedrep}, it is elementary to see that $M_{**} \cong K_*(\EMH_*)$ as $K_*$-modules.

The Hopf ring $K_*(\EMH_*)$ is generated as an $K_*$-algebra by 
primitives and one group-like element $a^{\Delta[0]} =_\text{def} [1] - [0]_0$. To prove the claims about the coalgebra structure, we only have to notice that $a^I$ is primitive when $i_0=1$. Clearly, $a_{(0)}$ is primitive. If $x$ is primitive and $y$ is any other element in the augmentation ideal, then
\begin{align*}
\Psi(x \circ y) = \Psi(x) \circ \Psi(y) & = \sum_{(y)} \left(x \circ y' \otimes [0] \circ y'' + [0] \circ y' \otimes x \circ y''\right) \\
&= x \circ y \otimes [0] + [0] \otimes x \circ y,
\end{align*}
and thus $x \circ y$ is also primitive. This shows that all elements of the form $a^I$ with $i_0=1$ are primitive. Conversely, if $i_0 = 0$, then $\Psi^p(a^I) = \left(a^{sI}\right)^{\otimes 
p}$ (mod $*$), and thus $a^I$ is not primitive.
\end{proof}

\begin{corollary} \label{kncohomology}
Choose a basis of $K_*(\EMH_*)$ containing the generators $a^I$ in $M_{**}$.
In the dual basis, let $x_{I}$ be the dual of $a^I$. Denote by 
$\widetilde{\EMH}$ the connected component of $0$ in $\EMH$. Then modulo phantoms,
\begin{equation} \label{kncohomologyformula}
K^*(\widetilde{\EMH}_*) \cong \smashop\bigotimes_{i_0=1} P_{t_0(I)+1}\left(x_I \right).
\end{equation}
\end{corollary}
\begin{proof}
For a Hopf algebra $A$ with a chosen basis, denote by $P(A)$ its primitives and by $Q(A)$ 
its indecomposables as a submodule (using the basis).

First note that $P(M_{**}) \subset Q(M_{**})$. (This is not true in $K_*(\EMH_*)$!)

Now $K^*(\widetilde{\EMH}_*)$ is precisely the sub-co-Hopf ring generated by the indecomposables, i.~e., the duals of the $*$-primitives. 

Thus the classes $x_I$ with $i_0 =1$ generate $K^*(\widetilde{\EMH}_*)$. The algebra structure follows by 
inspection of the coproduct.
\end{proof}

In any commutative and cocommutative Hopf algebra over $\F_p$ with suitable finiteness hypothesis, there are Frobenius and Verschiebung maps corresponding to the $p$-fold product and coproduct. In the case of a Hopf ring, they interact with the circle product in a simple way \cite[Section 7]{ravenel-wilson:morava}. We recall the definitions and basic properties for the reader's convenience.

\begin{defn}
For a graded algebra $A$ over $\F_p$, define the \emph{Frobenius} homomorphism $F\colon A \to A$ to be the $p$th power map $x \to x^{(1)}$. Similarly, for a graded coalgebra $C$ over $\F_p$ which, as a graded vector space, is the colimit of coalgebras $C_n$ of finite type, define the \emph{Verschiebung} $V\colon C \to C$ to be the continuous dual of the Frobenius on the pro-finite type algebra $C^{\vee}$.
\end{defn}

For H-spaces $X$ of finite type, the filtration $C_n$ of $K_*(X)$ is by definition given by the skeletal filtration of $X$. The spaces $\EMH_k$ are of finite type.

\begin{lemma}[{\cite[Lemma 7.1]{ravenel-wilson:morava}}] \label{frobeniusandverschiebung}
In a commutative and cocommutative Hopf ring $A$,
\begin{enumerate}
\item $V$ and $F$ are Hopf algebra maps multiplying resp. dividing the degree by $p$;
\item $VF(x)=FV(x)=[p]\circ x$;
\item $V(x \circ y) = V(x) \circ V(y)$;
\item $\Psi^p(x) \equiv V(x) \otimes \cdots \otimes V(x) \pmod{\text{asymmetric terms}}$;
\item $F(V(x) \circ y) = x \circ F(y)$. \label{FVcondition}
\end{enumerate}
\end{lemma}

The following results are dual to the pairing of bar spectral sequences introduced in \cite{thomason-wilson:pairing} or, in the case of the Eilenberg-Mac~Lane spectrum, \cite[Section 1]{ravenel-wilson:morava}.

\begin{prop}[Module structures on EMSS]\label{emssmodule}
Given a map of spaces $X \times T \xrightarrow{\mu} Y$ and a field spectrum $K$. Let $E^n(X)$ denote the $K$-based EMSS for the path-loop fibration on $X$, similarly for $Y$. Then there are homomorphisms
\[
E^n(X) \otimes_{K_*} K_*(T) \xrightarrow{\tilde\mu} E^n(Y)
\]
such that
\[
d^n(\tilde\mu(x,\eta)) = \tilde\mu(d^n(x),\eta) \quad \text{for $x \in E^n(X)$, $\eta \in K_*(T)$},
\]
and if $x$ is a permanent cycle representing a class $\xi \in K_*(\Omega X)$ then $\tilde\mu(x,\eta)$ is also permanent cycle and represents the class $\Omega\mu_*(\xi,\eta) \in \Omega Y$.

In the cobar resolution, $\tilde\mu$ is given by
\[
[x_1 \mid \cdots \mid x_s] \otimes \eta \mapsto \sum_{(\eta)} [\mu_*x_1,\eta(1)) \mid\cdots \mid \mu_*(x_s,\eta(s))],
\]
where the sum is given by $\Psi^s(\eta) = \sum_{(\eta)} \eta(1) \otimes \cdots \otimes \eta(s)$.
\end{prop}

\begin{proof}
Denote by $X^\bullet$ (resp. $Y^\bullet$) the cobar construction $C_X(*,*)$ (resp. $C_Y(*,*)$). Then the canonical map
\[
X^s \times T \xrightarrow{\id \times \diag^s} X^s \times T^s \to (X \times T)^s \xrightarrow{\mu^s} Y^s
\]
induces a map of cosimplicial spaces. We thus obtain a map of total towers
\[
\left(\Tot^s X^\bullet\right) \times T \to \Tot^s \left( X^\bullet \times T\right) \to \Tot^s Y^\bullet
\]
which, on homotopy inverse limits, agrees with the standard map
\[
(\Omega X) \times T \xrightarrow{\Omega\mu} \Omega Y \quad \text{given by} \quad (\Omega\mu)(\gamma,\eta)(t) = \mu(\gamma(t),\eta).
\]

Applying $K[-]$, we obtain the following diagram:
\[
\xymatrix{
K[(\Tot^sX^\bullet) \times T] \ar[r] \ar[d] & K[\Tot^s Y^\bullet] \ar[d]\\
\Tot^sK[X^\bullet \times T] \ar[r] & \Tot^sK[Y^\bullet].
}
\]
Since
\[
\pi_*\Tot^sK[X^\bullet \times T] \cong \pi_*\Tot^sK[X^\bullet] \otimes_{K_*} K_*(T),
\]
we obtain a commutative square
\[
\xymatrix{
K_*(\Tot^s X^\bullet) \otimes_{K_*} K_*(T) \ar[r] \ar[d] & K_*(\Tot^s Y^\bullet) \ar[d]\\
D^1_{s*}(X) \otimes_{K_*} K_*(T) \ar[r] & D^1_{s*}(Y)
}
\]
and hence a map of spectral sequences compatible with the filtration. The description of $\tilde \mu$ in the cobar complex follows directly from this construction, as the $s$-fold diagonal on $X$ is used in defining $\tilde \mu$. 
\end{proof}

%

Let $K$ be a commutative ring spectrum, and denote by ${}^rE^n_{s,t}$ the $E$-based EMSS for the path-loop fibration on $\EMH_r$; we write $E^*_{**}$ for the collection of all such spectral sequences for $r \geq 0$.

\begin{corollary} \label{pairing}
There is a homomorphism
\[
\circ\colon{}^rE^m_{s,t} \otimes K_{t'}(\EMH_{r'}) \to {}^{r+r'}E^m_{s,t+t'}
\]
compatible with the circle product
\[
\circ\colon K_{t-s}\EMH_r \otimes K_{t'}\EMH_{r'} \to K_{t+t'-s}\EMH_{r+r'}
\]
For $x \in E^n$ and $\eta \in K_*\EMH_*$,
\[
d^n(x \circ \eta) = d^n(x) \circ \eta.
\]
If $x$ is a permanent cycle representing a class $\xi \in K_*\EMH_*$, then $x \circ y$ is also a permanent cycle, and it represents $\xi \circ y$.

\end{corollary}
\begin{proof}
This is a special case of Prop.~\ref{emssmodule} for $C=\EMH_r$, $X=\EMH_r'$, $D=\EMH_{r+r'}$, and $\mu = \circ$.
\end{proof}

The following proposition follows from the multiplicative pairing and the relation between Frobenius and Verschiebung stated in Lemma~\ref{frobeniusandverschiebung}\eqref{FVcondition}:

\begin{prop} \label{filtrationjumpprop}
Let $x \in {}^rE^n_{**}$ and $\eta \in K_*(\EMH_{r'})$. Assume that $d^n(x) = z^{(1)}$ for some permanent cycle $z \in E^n_{**}$ representing a class $\zeta \in K_*(\EMH_r)$. Then there is an $m \geq n$ and a permanent cycle representing $\zeta \circ V(\eta)$ such that $x \circ \eta$ is an $(m-1)$-cycle and
\[
d^m(x \circ \eta) = t^{(1)}.
\]
\end{prop}

\begin{prop}[Determination of the EMSS for $\EMH_*$] \label{pnspectralsequence}
We have:
\begin{enumerate}
\item \label{torterm} As $K_*$-Hopf algebras,
\[
\Tor_{**}(K^*\widetilde{\EMH}_*) \cong 
\smashop \bigotimes_{i_0=1} \left(
\exterior\left(\sigma x_I \right) \otimes  \Gamma \left(\phi x_I \right)\right).
\]
\item \label{cotorterm} Dually, $E_{**}^2 = 
\Cotor_{**}(K_*\EMH_*)$ is given by
\[
E^2_{**} \cong  \bigotimes_{i_0=1} \left(\exterior\left(\sigma a^I \right) \otimes  \left(\phi a^I \right)\right).
\]
\item \label{transpotencecomp}
The cotranspotence $\phi a^I$ is represented in the cobar complex by the class
\[
\phi a^I \equiv \sum_{l=1}^{p-1} u_l \biggl[ 
\left(a^{s^{-t_0(I)}I}\right)^{*l} \biggm| 
\left(a^{s^{-t_0(I)}I} \right)^{*(p-l)} \biggr]
\]
for units $u_l \in \F_p^\times$, modulo classes that are more than $p$ times decomposable.

Furthermore, setting $m=l_0(I)$, we have that
\[
\phi (a_{(0)}) \circ a^I = \begin{cases}
0; & \text{if } m=0\\
\phi \left(a^{s^m(I) + \Delta_{n-m}}\right); & \text{otherwise}
\end{cases}
\]
or, equivalently,
\[
\phi (a^I) = \phi(a_{(0)}) \circ \left( 
a^{s^{-(t_0(I)+1)}I} \right)
\]
\item \label{differentials} The only sources of differentials are the 
factors of the form $\sigma a^I$ with $i_0=1$:
\begin{enumerate}
\item \label{difffirstkind} If $l=\Delta[m]$ for some $m\geq 1$ then
\[
d^{2p^m-1}(\sigma a^I) = \phi(a^I)^{(m)}
\]
\item \label{diffsecondkind} Otherwise, let $m = l_1(I)$ be the number of leading ones in $I$ and $m' = l_0(s^mI)$ be the number of zeroes following the $m$ ones in $I$. Then
\[
d^{2p^m-1}(\sigma a^I) = \left(\phi a^{c^{m+k-1}sI + \Delta_{n-m'}}\right)^{(m)}
\]
\end{enumerate}
\item \label{einftyterm} The spectral sequence collapses at $E^{2p^n}$ with 
\[
E^{2p^n}_{**} = E^\infty_{**} = \bigotimes_{i_0=1} P_{t_1(s^{-t_0(I)}I)} \phi(a^I).
\]
\item \label{einftyrepr} In $E_{**}^\infty$,
 $a^I$ is represented by $\phi a^{s^{l_0(I)I}+\Delta_{n-l_0(I)}}$ if $i_0=0$.
\end{enumerate}
\end{prop}

\begin{proof}

\noindent (\ref{torterm} and \ref{cotorterm}): This is a routine calculation using 
Corollary \ref{kncohomology}, using the basic building blocks from Lemma~\ref{basictor}.

\noindent \eqref{transpotencecomp}: If $R=P_k(x)$ is divided power algebra on $x$, then the transpotence element in the cobar complex is given by any one of the homologous representatives $[ x^l \mid x^{p^k-l} ]$. Dually, the cotranspotence is given by the sum of the duals of these classes. Thus
\begin{equation} \label{transpotencecobarrep}
\phi(a^I) = \sum_{k=1}^{p^{t_0(I)+1}-1} \biggl[ \left( x_I^k 
\right)^\vee \biggm| \left( x_I^{p^{t_0(I)}-k} \right)^\vee 
\biggr].
\end{equation}

Now if $k=k_0 + k_1 p + \cdots + k_r p^r$ with $0 \leq k_i <p$, the 
dual of $x_I^k$ is decomposable as
\[
\left(x_I^k\right)^\vee = \frac1{k_0!k_1!\cdots k_r!} 
\left(a^I\right)^{*k_0} * 
\left(a^{s^{-1}I}\right)^{*k_1} * \cdots
* \left(a^{s^{-r}I}\right)^{*k_r}.
\]
If we denote by $Q(k)=k_0 + \cdots + k_r$ the sum of the $p$-digits 
of $k$, the dual of $x_I^k$ is thus decomposable into $Q(k)$ 
factors. In order for the $k$-summand in 
\eqref{transpotencecobarrep} to have $p$ or less factors, we 
therefore need that $Q(k)+Q(p^{t_0(I)+1}-k) \leq p$. This is the case 
if and only if $k=l \, p^{t_0(I)}$ for $1 \leq l \leq p-1$, proving 
the first formula of part \eqref{transpotencecomp}.

For the second formula, first note that $\phi(a_{(0)}) = 
\sum_{l=1}^{p-1} u_l \bigl[ a_{(n-1)}^{*l} \bigm| a_{(n-1)}^{*(p-l)} \bigr]$ 
modulo $(p+1)$-decomposables. Now for $I$ with $i_0=0$,
\[
\bigl[a_{(n-1)}^{*l} \bigm| a_{(n-1)}^{*(p-l)} \bigr] \circ a^I = 
\biggl[ \left(a_{(n-1)} \circ a^{sI}\right)^{*l} \biggm| \left( 
a_{(n-1)} \circ a^{sI}\right)^{*(n-l)} \biggr]
\]
The sum of the right hand sides can be written (up to $(p+1)$-decomposables) as $\phi(a^{I'})$, where $I'$ is the multi-index $sI+\Delta_{n-1}$, shifted to the left such that $i'_0=1$. This amount is given by the number $m-1$ in the statement. The last formula of \eqref{transpotencecomp} follows by reindexing.

\noindent \eqref{einftyrepr}: By the computation of the EMSS for $\EMH_1$, $\phi(a_0)$ represents the class $[1]-[0] = a^{\Delta[0]} \in K_0(\EMH_0)$. Thus $a^I = ([1]-[0]) \circ a^Ib^J$ is represented by
\[
(\phi a_0) \circ a^I  = \phi a^{s^{l_0(I)}I + \Delta_{n-l_0(I)}} \quad \text{for $i_0 =0 $.}
\]

\noindent \eqref{differentials}: We have computed that $d^{2p-1}(\sigma a_{(0)}) = (\phi a_{(0)})^{(1)}$. By \eqref{einftyrepr}, the class $\phi a_{(0)}$ represents $[1]-[0]$.

We are thus in the situation of Prop.~\ref{filtrationjumpprop}, and there is a $j$ and a permanent cycle $t$ representing $([1]-[0]) \circ a^{sI}= a^{sI}$ such that
\[
d^j(\sigma a^I ) = d^j(\sigma a_{(0)} \circ a^{I-\Delta_0}) = t^{(1)}.
\]
We now apply Lemma \ref{normalizedrep} to $a^{sI}$ to get:
\[
a^{sI} = \pm \left(a^{c^{m-1}sI}\right)^{(m-1)},
\]
where $m=l_1(I)$. We will study which class represents $a^{c^{m-1}sI}$.

Write 
\[
I = \overset m {\overbrace{1\cdots1}}\overset {m'} {\overbrace{0\cdots0}}I',
\]
where $m \geq 1$, $m' \geq 0$, and $I'$ is either empty or starts with a $1$, and if $m'=0$ then $I$ is empty. Then
\[
c^{m-1}sI = \overset {m'} {\overbrace{0\cdots0}}I' 0 \overset {m-1} {\overbrace{1\cdots1}},
\]
and thus \eqref{einftyrepr} tells us that $a^{c^{m-1}sI}$ is represented by $\phi(a^J)$ with
\[
J = \begin{cases}
I' 0 \overset {m} {\overbrace{1\cdots1}} \overset {m'-1} {\overbrace{0\cdots0}}; & I' \text{ nonempty}\\
\overset {m} {\overbrace{1\cdots1}} \overset {m'} {\overbrace{0\cdots0}}; & I' \text{ empty}.
\end{cases}
\]
This determines all the differentials. For \eqref{einftyterm}, note that any $I$ with $i_0=1$ is of the form $J$ as above, with uniquely determined numbers $m$, $m'$ and subindices $I'$. This means that all the classes $\phi(a^I)$ are torsion (of order $p^{t_1(s^{-t_0(I)}I)}$) in $E^\infty$ and that all the classes $\sigma(a^I)$ support differentials.
\end{proof}

\begin{thm} \label{cohopfconvergence}
Let $n > 0$ and $t \neq n+1$ be integers. Then the Eilenberg-Moore spectral sequence
\[
E^2_{r,s} = \Cotor^{K(n)_*\left(\widetilde{\EMH}_t\right)}_{r,s} 
(K(n)_*,K(n)_*) \Longrightarrow K(n)_*\left(\EMH_{t-1}\right)
\]
converges $2p^n$-constantly.
\end{thm}
\begin{proof}
From Prop.~\ref{pnspectralsequence}\eqref{einftyterm} we know that
\[
E^{2p^n}_{**} = E^\infty_{**} = \bigotimes_{i_0=1} P_{t_1(s^{-t_0(I)}I)} \phi(a^I).
\]
and from Ravenel-Wilson's computation (Theorem~\ref{wilsonthm}), we know that
\[
K_*(\EMH_*) \cong \bigotimes_{i_0 = 0} P_{t_1(I)+1}(a^I) \otimes R(a^{\Delta[n]}),
\]
Furthermore, we know by Prop.~\ref{pnspectralsequence}\eqref{einftyrepr} that the comparison map from the latter to the target of the spectral sequence is such that $a^I$ (for $I \neq \Delta[n]$) is represented by $\phi a^{s^{l_0(I)I}+\Delta_{n-l_0(I)}}$. Since these classes have the same multiplicative order, we observe that the comparison map is an isomorphism.

Note that the spectral sequence does not converge for $t=n+1$ (the class $a^{\Delta[n]} \in K_*(\EMH_n)$ has no representative in $E^\infty$), which is not surprising since $K_*(\EMH_{n+1}) = 0$. However, for $t>n+1$, the EMSS again converges for trivial reasons: source and target are trivial.
\end{proof}

\begin{proof}[Proof of Theorem~\ref{hopfringconvergence}]
In Theorem~\ref{cohopfconvergence} we proved pro-constant convergence, and since $K(n)_*(\EMH_r)$ is a finite $K(n)_*$-module for all $r$ by the calculations of Ravenel and Wilson, $\{K(n)_*(X)\}_{X \in \Fin(\EMH_r)}$ is ind-constant. Thus we have Ind-convergence.
\end{proof}


\begin{thebibliography}{Tam94}

\bibitem[AM69]{artin-mazur:etale}
M.~Artin and B.~Mazur.
\newblock {\em Etale homotopy}.
\newblock Lecture Notes in Mathematics, No. 100. Springer-Verlag, Berlin, 1969.

\bibitem[BK72]{yellowmonster}
A.~K. Bousfield and D.~M. Kan.
\newblock {\em Homotopy limits, completions and localizations}.
\newblock Springer-Verlag, Berlin, 1972.
\newblock Lecture Notes in Mathematics, Vol. 304.

\bibitem[Boa99]{boardman:conditionally}
J.~Michael Boardman.
\newblock Conditionally convergent spectral sequences.
\newblock In {\em Homotopy invariant algebraic structures (Baltimore, MD,
  1998)}, volume 239 of {\em Contemp. Math.}, pages 49--84. Amer. Math. Soc.,
  Providence, RI, 1999.

\bibitem[Bou87]{bousfield:homology-cosimplicial}
A.~K. Bousfield.
\newblock On the homology spectral sequence of a cosimplicial space.
\newblock {\em Amer. J. Math.}, 109(2):361--394, 1987.

\bibitem[CE99]{ce:homalg}
Henri Cartan and Samuel Eilenberg.
\newblock {\em Homological algebra}.
\newblock Princeton Landmarks in Mathematics. Princeton University Press,
  Princeton, NJ, 1999.
\newblock With an appendix by David A. Buchsbaum, Reprint of the 1956 original.

\bibitem[Dol71]{dold:chern-classes}
Albrecht Dold.
\newblock Chern classes in general cohomology.
\newblock In {\em Symposia Mathematica, Vol. V (INDAM, Rome, 1969/70)}, pages
  385--410. Academic Press, London, 1971.

\bibitem[Dwy74]{dwyer:strongconvergence}
W.~G. Dwyer.
\newblock Strong convergence of the {E}ilenberg-{M}oore spectral sequence.
\newblock {\em Topology}, 13:255--265, 1974.

\bibitem[Dwy75]{dwyer:exoticconvergence}
William~G. Dwyer.
\newblock Exotic convergence of the {E}ilenberg-{M}oore spectral sequence.
\newblock {\em Illinois J. Math.}, 19(4):607--617, 1975.

\bibitem[EM66]{eilenberg-moore:homology-and-fibrations}
Samuel Eilenberg and John~C. Moore.
\newblock Homology and fibrations. {I}. {C}oalgebras, cotensor product and its
  derived functors.
\newblock {\em Comment. Math. Helv.}, 40:199--236, 1966.

\bibitem[GJ99]{goerss-jardine}
Paul~G. Goerss and John~F. Jardine.
\newblock {\em Simplicial homotopy theory}.
\newblock Birkh\"auser Verlag, Basel, 1999.

\bibitem[Hod75]{hodgkin:emss}
Luke Hodgkin.
\newblock The equivariant {K}\"unneth theorem in {$K$}-theorem.
\newblock In {\em Topics in $K$-theory. Two independent contributions}, pages
  1--101. Lecture Notes in Math., Vol. 496. Springer, Berlin, 1975.

\bibitem[Isa04]{isaksen:strictmodel}
Daniel~C. Isaksen.
\newblock Strict model structures for pro-categories.
\newblock In {\em Categorical decomposition techniques in algebraic topology
  (Isle of Skye, 2001)}, volume 215 of {\em Progr. Math.}, pages 179--198.
  Birkh\"auser, Basel, 2004.

\bibitem[JO99]{jeanneret-osse:emss}
A.~Jeanneret and A.~Osse.
\newblock The {E}ilenberg-{M}oore spectral sequence in {$K$}-theory.
\newblock {\em Topology}, 38(5):1049--1073, 1999.

\bibitem[MS06]{may-sigurdsson}
J.~P. May and J.~Sigurdsson.
\newblock {\em Parametrized homotopy theory}, volume 132 of {\em Mathematical
  Surveys and Monographs}.
\newblock American Mathematical Society, Providence, RI, 2006.

\bibitem[RS65]{rothenberg-steenrod}
M.~Rothenberg and N.~E. Steenrod.
\newblock The cohomology of classifying spaces of {$H$}-spaces.
\newblock {\em Bull. Amer. Math. Soc.}, 71:872--875, 1965.

\bibitem[RW80]{ravenel-wilson:morava}
Douglas~C. Ravenel and W.~Stephen Wilson.
\newblock The {M}orava {$K$}-theories of {E}ilenberg-{M}ac {L}ane spaces and
  the {C}onner-{F}loyd conjecture.
\newblock {\em Amer. J. Math.}, 102(4):691--748, 1980.

\bibitem[Sey78]{seymour:emss}
R.~M. Seymour.
\newblock On the convergence of the {E}ilenberg-{M}oore spectral sequence.
\newblock {\em Proc. London Math. Soc. (3)}, 36(1):141--162, 1978.

\bibitem[Shi96]{shipley:homology-cosimplicial}
Brooke~E. Shipley.
\newblock Convergence of the homology spectral sequence of a cosimplicial
  space.
\newblock {\em Amer. J. Math.}, 118(1):179--207, 1996.

\bibitem[Smi70]{smith:emss}
Larry Smith.
\newblock {\em Lectures on the {E}ilenberg-{M}oore spectral sequence}.
\newblock Lecture Notes in Mathematics, Vol. 134. Springer-Verlag, Berlin,
  1970.

\bibitem[Tam94]{tamaki:emss}
Dai Tamaki.
\newblock A dual {R}othenberg-{S}teenrod spectral sequence.
\newblock {\em Topology}, 33(4):631--662, 1994.

\bibitem[TW80]{thomason-wilson:pairing}
Robert~W. Thomason and W.~Stephen Wilson.
\newblock Hopf rings in the bar spectral sequence.
\newblock {\em Quart. J. Math. Oxford Ser. (2)}, 31(124):507--511, 1980.

\bibitem[Wil84]{wilson:hopf-ring-morava}
W.~Stephen Wilson.
\newblock The {H}opf ring for {M}orava {$K$}-theory.
\newblock {\em Publ. Res. Inst. Math. Sci.}, 20(5):1025--1036, 1984.

\bibitem[Wil00]{wilson:hopf-rings}
W.~Stephen Wilson.
\newblock Hopf rings in algebraic topology.
\newblock {\em Expo. Math.}, 18(5):369--388, 2000.

\end{thebibliography}

\end{document}